# A Convolutional Hierarchical Deep-learning Neural Network (C-HiDeNN) Framework for Non-linear Finite Element Analysis

Yingjian Liu[1], Monish Yadav Pabbala[1], Jiachen Guo[3], Chanwook Park[2], Gino Domel[2], Wing Kam Liu[2], Dong Qian[1,*]

[1] Department of Mechanical Engineering, The University of Texas at Dallas, Richardson, TX 75080, USA

[2] Department of Mechanical Engineering, Northwestern University, 2145 Sheridan Rd., Evanston IL 60208, USA

[3] Theoretical and Applied Mechanics Program, Northwestern University, 2145 Sheridan Road, Evanston, 60201, IL, USA

*Corresponding author, e-mail: dong.qian@utdallas.edu

## Abstract

We present a framework for the Convolutional Hierarchical Deep-learning Neural Network (C-HiDeNN) tailored for nonlinear finite element analysis. Building upon the structured foundation of HiDeNN, C-HiDeNN introduces a convolution operator to enhance numerical approximation. A distinctive feature of C-HiDeNN is its higher-order accurate approximation achieved through an expanded set of parameters, such as the polynomial order '*p*,' dilation parameter '*a*,' patch size '*s*,' and nodal position 'X'. These parameters function as the functional equivalents of weights and biases within each C-HiDeNN patch. In addition, C-HiDeNN can be selectively applied to regions requiring high resolution to adaptively improve local prediction accuracy. To demonstrate the effectiveness of this framework, we provide numerical examples in the context of nonlinear finite element analysis. The results show that our approach achieves significantly higher accuracy than conventional Finite Element Method (FEM) while substantially reducing computational costs.

# 1 Introduction

Deep learning and neural networks [1, 2] methods have received a great deal of attention as a tool for enhancing computational speed. Its applications in different areas have been reported, e.g. image analysis [3-5], material design and discovery [6-8]. With the universal approximation (**UA**) capabilities, deep learning can replicate linear or non-linear functions [9, 10] and offers a new pathway for solving ordinary and partial differential equations (ODEs/PDEs) in science and engineering. Inspired by the UA capabilities, the Hierarchical deep-learning neural networks (HiDeNN) was proposed by Zhang et al 2021 [11]. The basic idea of HiDeNN is to establish building blocks of deep neural network (DNN) to replicate the shape function approximations based on both standard finite element (FE), meshfree methods, and others, such as those based on isogeometric analysis (IGA). In the case of standard FE, the resulting method is called HiDeNN-Finite Element Method. Once the approximation is established, the DNN can be further extended to incorporate a loss function that is directly related to the problem's physics. For instance, a strain energy potential is introduced as the loss function in the case of linear elasticity. Compared to regular FEM, an important feature of HiDeNN-FEM is that the weight and bias of the constructed DNN are also functions of the nodal coordinates. As a result, the training process is directly integrated with the solution process, and it naturally leads to optimization of the nodal coordinates. This is known as the *r*-adaptivity. The constructed DNN is also capable of adding additional neurons to refine the mesh. As such, HiDeNN-FEM is capable of performing *rh*-adaptivity. Application of HiDeNN-FEM for both linear and nonlinear problems has been demonstrated in [11-13], with extensions to tensor decomposition [14], reduced order model [15, 16], topology optimization [17], and convolution finite element method [18, 19]. In the case of nonlinear FE, it has been shown that HiDeNN-FEM effectively captures stress concentration, removes the volumetric locking and hourglass mode, and can handle extremely large deformation. All of these capabilities were achieved without the need for sophisticated treatment as being done in the case of regular FE.

A major step beyond HiDeNN-FEM is the introduction of convolution operators into the FE-based approximation. The convolution operator is motivated by convolutional neural networks (CNNs) [20] that have been extensively used for image analysis [21, 22]. As shown in the rest of this paper, this incorporation leads to a hybrid FE/meshfree approximation. The resulting method is called the Convolution Hierarchical Deep-learning Neural Network (C-HiDeNN) [18, 19]. C-HiDeNN interpolants exhibit two key properties: the Kronecker-delta property and the reproducing condition that allows the interpolant to approximate functions at arbitrary orders. The Kronecker-delta property simplifies the treatment of boundary conditions, which are typically more complex in meshfree methods. The reproducing condition enables the construction of higher-order, more accurate approximations without requiring additional degrees of freedom. The higher-order accurate approximation is achieved though introducing additional parameters include the polynomial order '*p*', dilation parameter '*a*', patch size '*s*', and mesh position '*x*'. The robustness of C-HiDeNN have been demonstrated extensively for transient heat transfer problems in additive manufacturing [14], and linear elasticity problems [11]. Park et al [18] demonstrated that the computational cost of C-HiDeNN can be significantly reduced with the use of graphic processing unit (GPU) and to the same order as the regular FEM. The solver developed also makes extensive use of the JAX library

((https://github.com/google/jax), which is a high-level Python library developed by Google and provides accelerated linear algebra library at the backend that can be implemented on CPU, GPU and TPU (tensor processing unit). Other significant extensions of C-HiDeNN include C-HiDeNN-Tensor Decomposition (TD) [14], which was aimed at developing reduced order models for high-dimensional problems. It was shown that C-HiDeNN-TD effectively converts an *n*-dimensional problem to *n* 1-dimensional problems, thus leading to computational cost scales linearly with the degree of freedom, as opposed to exponential growth. Applications of C-HiDeNN-TD for large-scale physical systems have been demonstrated in [23]. Comparison with regular solvers of partial differential equations showed significant gains in computing speed with much reduced memory consumption and data storage. On the other hand, this approach does not require the expensive offline data training process as the conventional data-driven approach.

Previous implementations of C-HiDeNN for solid mechanics have been limited to linear elasticity problems. This study presents a significant extension of C-HiDeNN to nonlinear finite element employing total Lagrangian formulation. As shown in the rest of the paper, the C-HiDeNN interpolant extends the limit of the conventional low-order FE approximation and contributes to solutions with much enhanced accuracy without the need for remeshing. Similar to the work of [18, 19], we leverage the high-performance capabilities of Google JAX (https://github.com/google/jax.) and JAX-FEM [24] method for automatic differentiation. This facilitates the computation of the material stiffness matrix and its application in solving inverse problems. To provide an understanding of how the interpolation function is constructed in C-HiDeNN, we first demonstrate examples of 1D by employing regular polynomial basis. We then extend to higher dimensions using

Radial Basis Function (RBF) and illustrate the application of C-HiDeNN with numerical examples.

The remainder of the paper is organized as follows: Section 2 reviews the basic concepts, formulation of C-HiDeNN and how neural network representations of C-HiDeNN are constructed. In section 3, algorithms of C-HiDeNN for nonlinear FE and the implementation are outlined. Section 4 provides several numerical examples to demonstrate the computational efficiency and generalization of the proposed method. Finally, conclusions and future directions are discussed in Section 5.

## 2 Nonlinear C-HiDeNN formulation

### 2.1 Non-linear FEM governing equations

Considering a solid initially occupying a bounded domain of $\Omega_0$ with the boundary of $\Gamma_0$, the governing equation of the conservation of linear momentum [25] is given as

$$\nabla_0 \cdot \mathbf{P} + \rho_0 \mathbf{b} = \rho_0 \ddot{\mathbf{u}}, \tag{1}$$

in which $\nabla_0$ is the "del" operator with respect to the material coordinates **X**, **P** is the first Piola-Kirchhoff stress, $\rho_0$ is the mass density, **b** is the body force, **u** is the displacement, and superimposed dot indicates the time derivative. The boundary $\Gamma_0$ consists of the essential boundary $\Gamma_{u0}$ on which the displacement is imposed and the natural boundary $\Gamma_{t0}$ on which the traction is applied.

The two most widely used approaches in nonlinear FEM for solids are the total Lagrangian (TL) formulation and updated Lagrangian formulation (UL). In this work, the total Lagrangian (TL) [25] is employed. C-HiDeNN can be applied in the context of UL in a similar way. In TL, the weak form of the momentum equation is expressed as:

$$\int_{\Omega_0} \rho_0 \ddot{\mathbf{u}} \cdot \delta \mathbf{u} \; d\Omega + \int_{\Omega_0} \mathbf{P} : \delta \mathbf{F}^{\mathrm{T}} d\Omega - \int_{\Omega_0} \mathbf{b} \cdot \delta \mathbf{u} \; d\Omega - \int_{\Gamma_{t0}} \mathbf{t} \cdot \delta \mathbf{u} \, d\Gamma = 0, \tag{2}$$

in which $\mathbf{u}$=$\mathbf{x}$-$\mathbf{X}$ denotes the displacement field with $\mathbf{X}$ being spatial coordinate. The symbol $\delta$ represents the variation, $\mathbf{F}$ represents the deformation gradient, given as $\mathbf{F}=\frac{\partial \mathbf{x}}{\partial \mathbf{X}}$. Additionally, $\mathbf{t}$ represents the traction applied on the boundary $\Gamma_{t0}$, i.e.,$\mathbf{t} = \bar{\mathbf{t}}_{\mathbf{0}}$ on $\Gamma_{\mathrm{t0}}$.

To solve the weak form Eq. (2), a Lagrangian mesh [25] is introduced, along with the C-HiDeNN shape function $\widetilde{N}_I$ defined at node indexed by $I$. A common approach is to introduce the parent domain and express $\widetilde{N}_I$ as a function of the parent element coordinates $\xi$, i.e., $\widetilde{N}_I = \widetilde{N}_I(\xi)$. The mapping from the parent to the initial and parent to the current configurations can then be constructed as $\mathbf{X} = \mathbf{X}(\xi)$ and $\mathbf{x} = \mathbf{x}(\xi, t)$, respectively. After substituting the approximation based on the FEM shape function into Eq. (2), the discretized momentum equation is given as

$$\mathbf{M}\ddot{\mathbf{u}} = \mathbf{f}^{ext} - \mathbf{f}^{int}, \tag{3}$$

with

$$\mathbf{f}^{int} = f_{iI}^{int} = \int_{\Omega_0} \frac{\partial \tilde{N}_I}{\partial X_j} P_{ji} d\Omega_0 = \int_{\Omega_0} \left( \mathcal{B}_{Ij}^0 \right)^T P_{ji} d\Omega_0, \tag{4}$$

$$\mathbf{f}^{ext} = f_{iI}^{ext} = \int_{\Omega_0} \tilde{N}_I \rho_0 b_i d\Omega_0 + \int_{\Gamma_{t_i}^0} \tilde{N}_I \bar{t}_i^{\,0} d\Gamma_0, \tag{5}$$

$$\mathbf{M} = M_{ijIJ} = \hat{\delta}_{ij} \int_{\Omega_0} \rho_0 \tilde{N}_I \tilde{N}_J d\Omega_0. \tag{6}$$

where lower case index indicates the dimension, $\mathbf{M}$ is the mass matrix defined by Eq.(6) in which $\hat{\delta}_{ij}$ is the Kronecker delta. The shape function derivative matrix is defined as $\mathcal{B}_{jI}^0 = \frac{\partial \tilde{N}_I}{\partial X_j}$.

## 2.2 C-HiDeNN interpolation

As can be seen from Eqs. (3) to (6), solving the linear momentum equation requires the evaluation of the shape function $\tilde{N}_I$ and its material derivatives $\frac{\partial \tilde{N}_I}{\partial X_j}$. In C-HiDeNN, a convolution kernel is introduced to incorporate information from neighboring elements in constructing the shape function and its derivatives (**Fig. 1**). As demonstrated later in this section, this results in a hybrid FE/meshfree representation, making C-HiDeNN a natural extension of the HiDeNN-FEM framework. In HiDeNN-FEM, the building blocks of neural networks (NN) consist of fundamental operations such as linear transformations, multiplication, and inversion (**Fig. 1**b), which are used to construct finite element (FE) shape functions. Since these operations are applied at the element level, they inherently lead to a trimmed neural network as represented by the blue block in **Fig. 1**a. C-HiDeNN further enhances the FE approximation by introducing a patch around an element (**Fig. 1**c) and applying a convolution operation to the FE shape function values generated from HiDeNN-FEM, enabling a more efficient and flexible representation. This operation is represented by the yellow block in **Fig. 1** (a). Based on the nature of the convolution operator, the neural network structure in the yellow block can be either trimmed, i.e., only neighboring elements are considered, or fully connected.

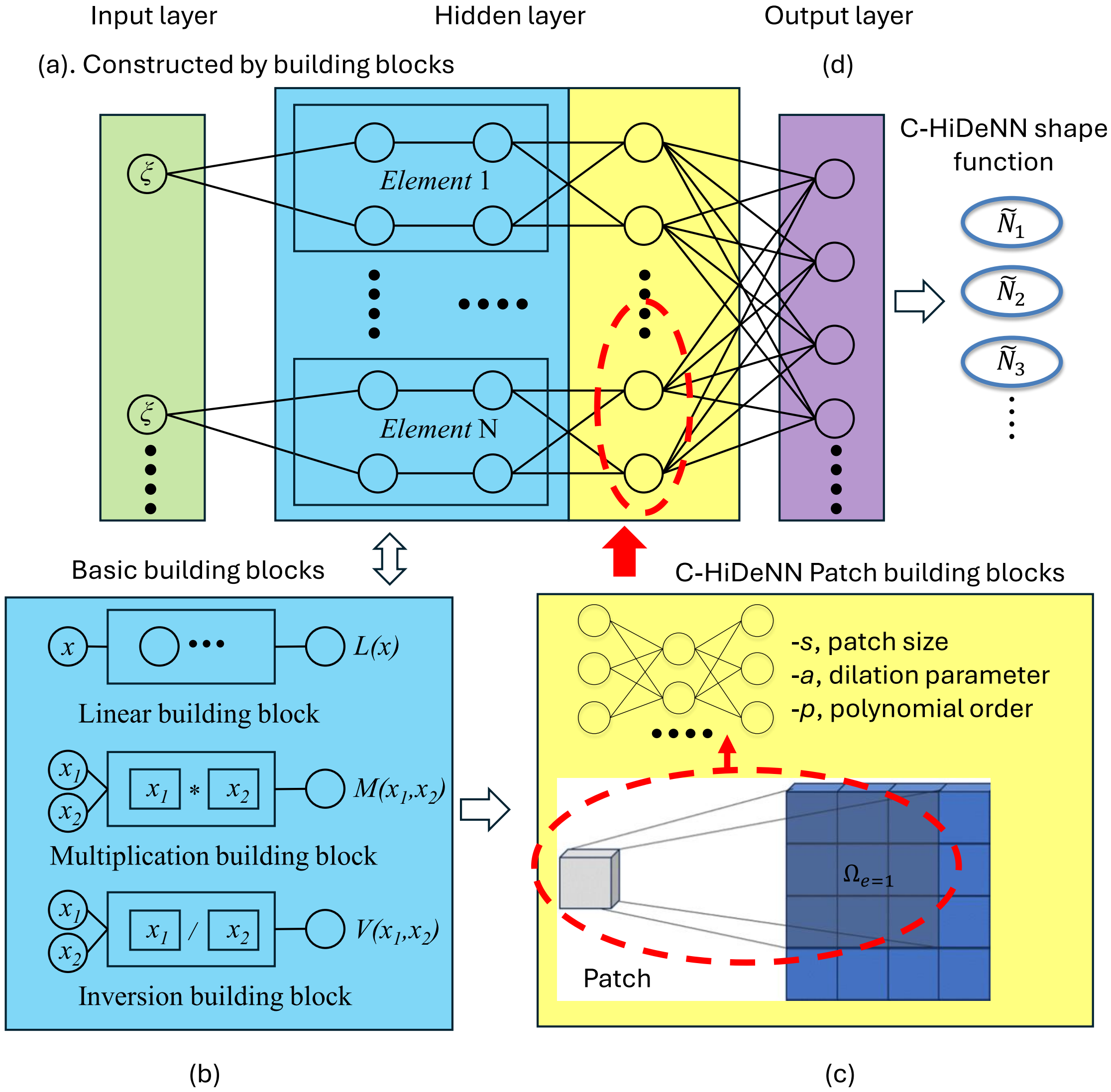

**Fig. 1. An illustration of the building blocks for constructing Convolution-HiDeNN patch function. (a) C-HiDeNN construction building blocks. (b) HiDeNN basic building blocks of linear, multiplication and inversion operations. (c) C-HiDeNN convolution patch. (d) C-HiDeNN output layer with a convolution kernel function.**

We start with a 1D illustration of the steps to construct C-HiDeNN interpolation. As shown in **Fig. 2**, a 1D parent domain is discretized by the standard 2-node linear elements with uniform nodal spacing of Δ and a set of $N$ nodes indexed by 1, 2 … $I$-2, $I$-1, $I$, $I$+1, $I$+2, …, $N$. The concept of the nodal patch $A_s^i$ is introduced and defined as

$$A_s^i = \left\{ J \in \mathbb{N} \middle| \left| \xi_J - \xi_i \right| \leq s \bullet \Delta \right\}. \tag{7}$$

In Eq. (7), $\xi_J$ and $\xi_i$ are the $J$-th and $i$-th element coordinates, respectively, parameter $s$ is called the patch size, which corresponds to the number of layers of elements surrounding node $i$. Next, we introduce the element patch, denoted by $A_s^e$, given as $A_s^e = \bigcup_{i \in A_e} A_s^i$ in which $A_e$ contains the set of nodes that belong to element $e$. Using the $k$-th element in **Fig. 2** as an example, and assuming $s$=1, we have $A_{s=1}^{i=I} = \{I-1, I, I+1\}$, $A_{s=1}^{i=I+1} = \{I, I+1, I+2\}$ and $A_{s=1}^{e=k} = \{I-1, I, I+1, I+2\}$.

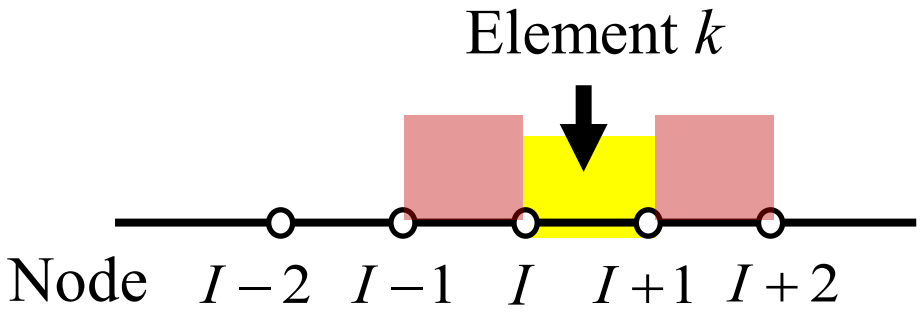

**Fig. 2. A 1D illustration of nodal and element patch.**

Based on the definitions given above, the field variable $u^h$ interpolated by Convolution HiDeNN-FEM is expressed as,

$$u^h(\xi) = \sum_{I \in A^e} N_I(\xi) \sum_{J \in A_s^I} \tilde{W}_{a,p,J}(X(\xi)) u_J = \sum_{K \in A_s^e} \tilde{N}_K(\xi) u_K, \tag{8}$$

in which $N_I(\xi)$ is the standard FE shape function, $\widetilde{W}_{a,p,J}$ is called a convolution patch function that is defined at node $J \in A_s^I$ and evaluated at the material coordinate $X$. The subscripts $a$ and $p$ indicate the dilation parameter and order of the approximation, respectively. The first part of the double summation in Eq. (8) resembles the FEM approximation, as it sums over the nodes within an element. The second part sums over the nodal patch for each node within the element. Since the nodal patch includes nodes outside of the element, the convolution patch function constructs a meshfree approximation. By using the element patch, the two summations can be combined, resulting in the last expression in Eq. (8).

An example of an application of Eq. (8) is shown in **Fig. 3**. We consider a 1D parent domain discretized by a set of uniformly spaced nodes. For simplicity we assume $X=\xi$ and employ linear element and quadratic interpolation for the convolution patch function, i.e., $a$=1, $p$=2. With the choice of $s$=1, we have $A_{s=1}^{I=1}=\{0,1,2\}$, $A_{s=1}^{I=2}=\{1,2,3\}$ and for element 1 we have $A_{s=1}^{e=1}=\{0,1,2,3\}$.

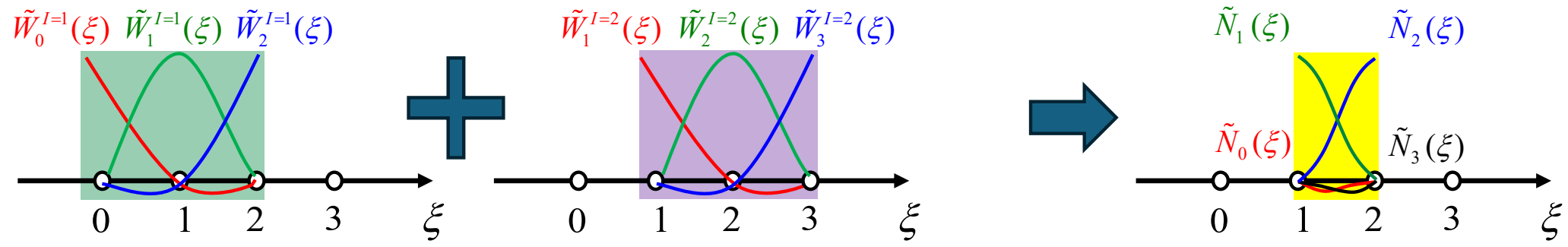

**Fig. 3. Construction of 1D C-HiDeNN shape function based on Lagrange polynomial kernel.**

For element 1, the C-HiDeNN shape functions are given below and plotted in **Fig. 3**:

$$
\begin{aligned}
\tilde{N}_0(\xi) &= -0.5(\xi-1)(\xi-2)^2, \\
\tilde{N}_1(\xi) &= \xi(\xi-2)^2+0.5(\xi-1)(\xi-2)(\xi-3), \\
\tilde{N}_2(\xi) &= -(\xi-1)^2(\xi-3)-0.5\xi(\xi-1)(\xi-2), \\
\tilde{N}_3(\xi) &= 0.5(\xi-1)^2(\xi-2).
\end{aligned} \tag{9}
$$

The convolution patch function is not limited to polynomial-based interpolations and can be extended by introducing more general meshfree types of interpolations. Examples include interpolations based on radial basis functions (RBF) [26], reproducing kernel particle method (RKPM) [27], isogeometric analysis [28] and others.

We further make an extension to convolution patch functions in multiple dimensions (2D and 3D) using radial basis functions (RBFs). Using the radial point interpolation approach [29], the patch function based on RBFs can be shown to be in the form of an $n$-component row vector $\tilde{\mathbf{W}}(\mathbf{X})$, taken from the $(n+m)$-component vector $\mathbf{W}(\mathbf{X})$, given by

$$\mathbf{W}(\mathbf{X})=\left[\mathbf{R}^T(\mathbf{X}) \quad \mathbf{P}^T(\mathbf{X})\right]\mathbf{G}^{-1} \tag{10}$$

in which $\mathbf{R}(\mathbf{X})$ is the matrix containing the radial basis functions evaluated over the nodes within a patch and is of size $n \times n$ with $n$ being the number nodes included in $A_s^I$, $\mathbf{P}(\mathbf{X})$ contains the polynomial basis function with order parameter $m$, e.g., $m=4$ in the case of 3D with basis functions of $1, x, y, z$. $\mathbf{G}$ is an assembled matrix and takes the form of

$$\mathbf{G}=\begin{bmatrix} \mathbf{R} & \mathbf{P} \\ \mathbf{P}^T & \mathbf{0} \end{bmatrix} \tag{11}$$

Based on Eq.(10), the derivative of the patch function can be obtained by

$$\mathbf{W}'(\mathbf{X})=\left[\left(\mathbf{R}^T(\mathbf{X})\right)' \quad \left(\mathbf{P}^T(\mathbf{X})\right)'\right]\mathbf{G}^{-1} \tag{12}$$

in which the superscript of prime indicates the differentiation with respect to a particular coordinate. For details on solving the convolution patch functions and derivatives, we refer to [29].

A case of generating the C-HiDeNN shape functions in 2D is shown in **Fig. 4** by employing the radial basis function as the patch function by setting the patch parameters $s$=1, $p$=1 and linear basis ( $m=3$ ). In **Fig. 4** the 2D domain is discretized by the standard 4-node quadrilateral elements with a uniform nodal spacing of $\Delta$. The nodal patch $A_s^I$ is defined as

$$A_s^I=\left\{J \in \mathbb{N} \,\middle|\, \left|\xi_J-\xi_I\right| \le s\bullet\Delta, \quad \left|\eta_J-\eta_I\right| \le s\bullet\Delta\right\}, \tag{13}$$

where, $\xi_J$, $\xi_I$, $\eta_J$, and $\eta_I$ are node $J$ and node $I$'s element coordinates in 2D. "$s$" is patch size.

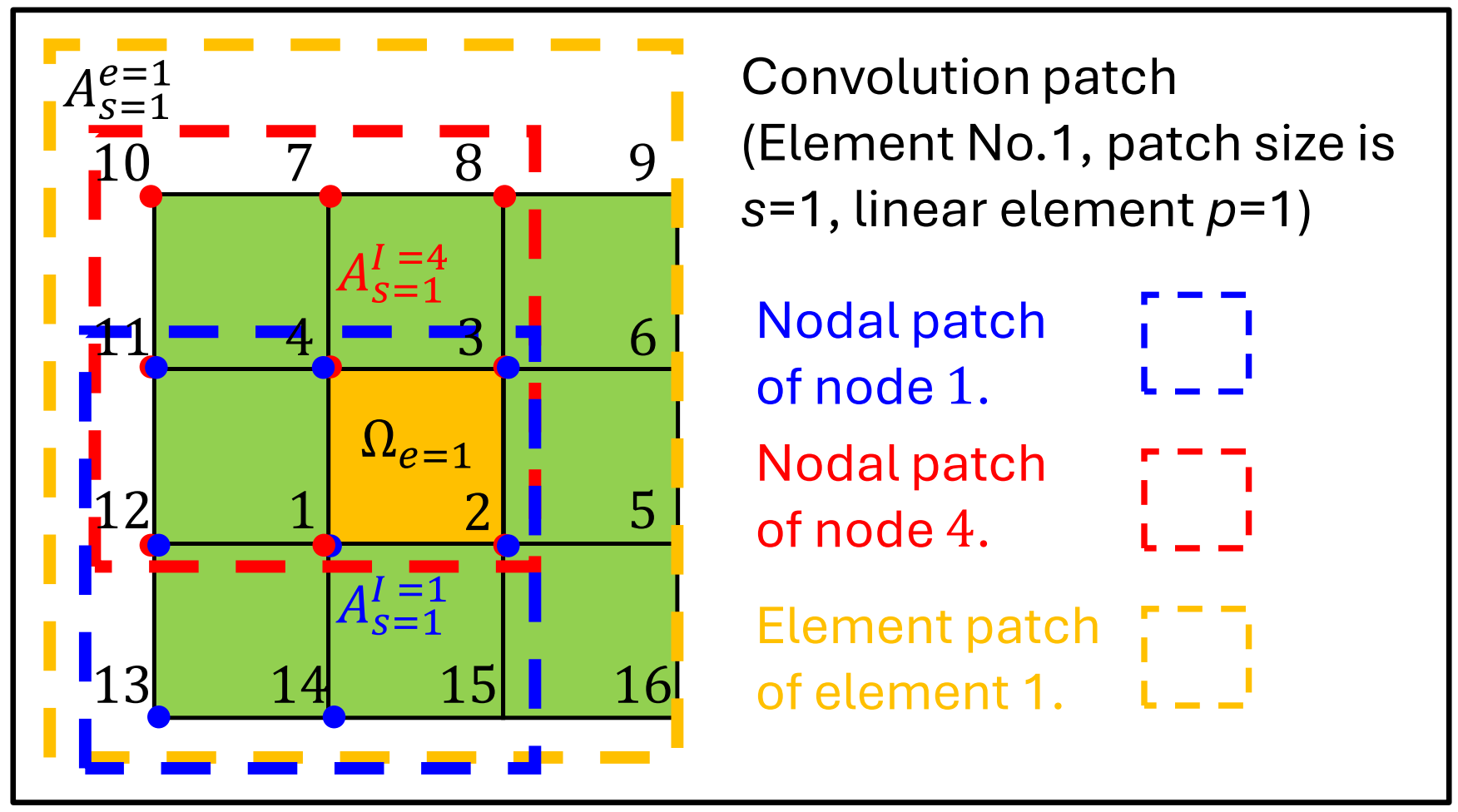

**Fig. 4. 2D C-HiDeNN patch illustration with patch size 1.**

A comparison between the standard linear FEM shape function and the C-HiDeNN shape function is illustrated in **Fig. 5**. In **Fig. 5**, a C-HiDeNN patch with a patch size of 1 and a linear polynomial basis function. **Fig. 5** (a) to **Fig. 5** (d) show the standard 4-node quadrilateral linear FE shape functions $N_i$, ranging from node 1 to node 4. **Fig. 5** (e) to **Fig. 5** (t) shows the corresponding C-HiDeNN shape functions $\widetilde{N}_i$, ranging from node 1 to node 16. The nodal indices are according to **Fig. 4**.

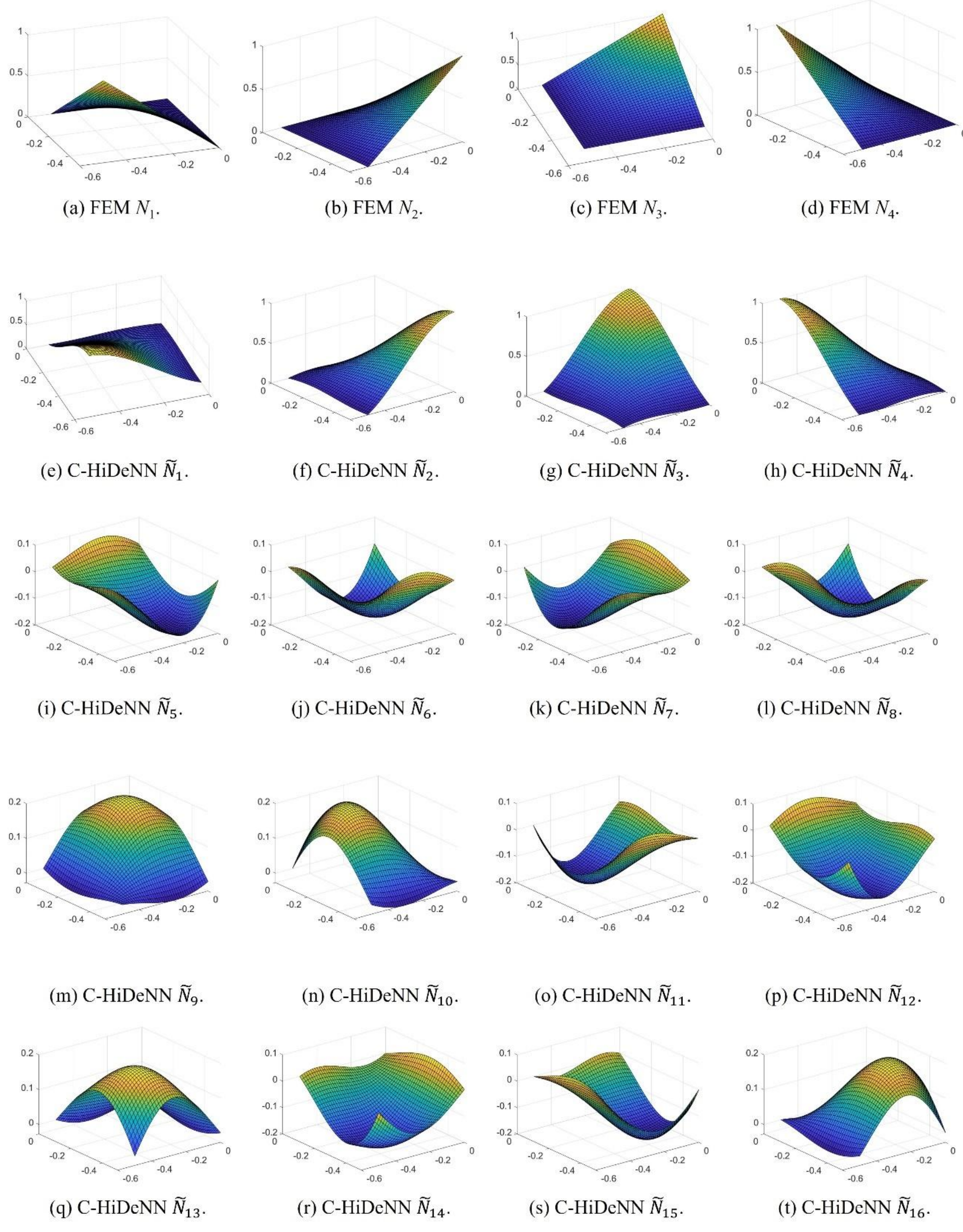

(a) FEM $N_1$. (b) FEM $N_2$. (c) FEM $N_3$. (d) FEM $N_4$.

(e) C-HiDeNN $\widetilde{N}_1$. (f) C-HiDeNN $\widetilde{N}_2$. (g) C-HiDeNN $\widetilde{N}_3$. (h) C-HiDeNN $\widetilde{N}_4$.

(i) C-HiDeNN $\widetilde{N}_5$. (j) C-HiDeNN $\widetilde{N}_6$. (k) C-HiDeNN $\widetilde{N}_7$. (l) C-HiDeNN $\widetilde{N}_8$.

(m) C-HiDeNN $\widetilde{N}_9$. (n) C-HiDeNN $\widetilde{N}_{10}$. (o) C-HiDeNN $\widetilde{N}_{11}$. (p) C-HiDeNN $\widetilde{N}_{12}$.

(q) C-HiDeNN $\widetilde{N}_{13}$. (r) C-HiDeNN $\widetilde{N}_{14}$. (s) C-HiDeNN $\widetilde{N}_{15}$. (t) C-HiDeNN $\widetilde{N}_{16}$.

**Fig. 5. 2D 4-node quadrilateral element comparison between FEM linear shape function and C-HiDeNN reconstructed higher order shape function with RBF kernel: (a-d) FEM shape function $N_1$, $N_2$, $N_3$ and $N_4$, (e-t) C-HiDeNN shape function $\widetilde{N}_1$ to $\widetilde{N}_{16}$. Nodes are indexed according to Fig.4.**

### 2.3 Neural Network Representation in C-HiDeNN

The process of constructing the neural network representation of the 1D C-HiDeNN shape functions is illustrated in **Fig. 6**. Here we consider a 1D domain with nodes indexed by …, -1,0,1,2, … We consider an element indexed by $e=1$ as shown in **Fig. 6**. For the neural network representation, the weight function and bias are defined as $W_{ij}^{lm}$ and $b_j^m$, where $i$ and $j$ denote the index of the neurons in layer $l$ and $m$, respectively. The activation function is given as, $\boldsymbol{\mathcal{A}}(x)=x$. Other activation function such as the ReLU function, i.e, $\boldsymbol{\mathcal{A}}_1(x)=\mathbf{max}(x,0)$, are used in the C-HiDeNN framework to construct other building blocks. For instance, the multiplication building block $M$ employs ReLU [12]. For a given element coordinate $\xi$, we first construct the linear FE shape function based on HiDeNN-FEM as shown in the top row of **Fig. 6**. For element 1, the corresponding NN representations for the regular linear shape functions are given as

$$
\begin{aligned}
N_1(\xi;\mathbf{W},\mathbf{b},\boldsymbol{\mathcal{A}}) &= W_{11}^{23}\boldsymbol{\mathcal{A}}\left(W_{11}^{12}\xi+b\right) \\
N_2(\xi;\mathbf{W},\mathbf{b},\boldsymbol{\mathcal{A}}) &= W_{22}^{23}\boldsymbol{\mathcal{A}}\left(W_{22}^{12}\xi+b\right)
\end{aligned} \tag{14}
$$

The specific values of the weight and bias are shown in **Fig. 6**.

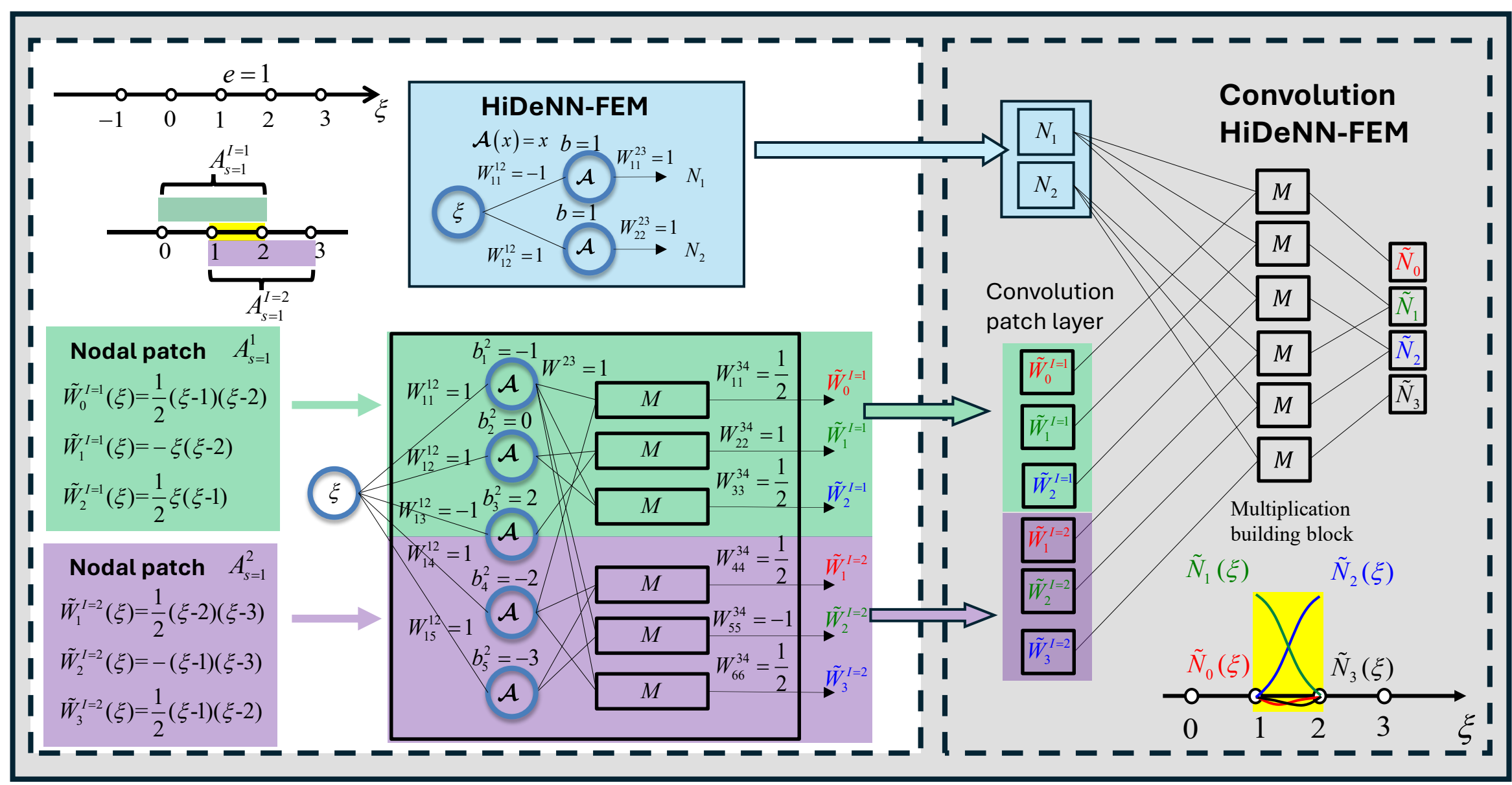

**Fig. 6. An illustration of the 1D C-HiDeNN formulation and neural network structure. The input coordinates are from two different element patches and serve as inputs to the HiDeNN-FEM building blocks. The convolution kernel performs the convolution operation to form the final C-HiDeNN shape function.**

The same element coordinate can also be employed as an input for generating the convolution patch functions as shown in the second row in **Fig. 6**. In the case of Lagrangian polynomial with $p$ = 2, the NN representations are expressed as

$$
\begin{aligned}
\tilde{W}_0^{I=1}\left(\xi;\mathbf{W},\mathbf{b},\mathcal{A}\right) &= W_{11}^{34}\left\{M\left[(W^{23}\mathcal{A}(W_{11}^{12}\xi+b_1^2)),(W^{23}\mathcal{A}(W_{14}^{12}\xi+b_4^2))\right]\right\},\\
\tilde{W}_1^{I=1}\left(\xi;\mathbf{W},\mathbf{b},\mathcal{A}\right) &= W_{22}^{34}\left\{M\left[(W^{23}\mathcal{A}(W_{12}^{12}\xi+b_2^2)),(W^{23}\mathcal{A}(W_{13}^{12}\xi+b_3^2))\right]\right\},\\
\tilde{W}_2^{I=1}\left(\xi;\mathbf{W},\mathbf{b},\mathcal{A}\right) &= W_{33}^{34}\left\{M\left[(W^{23}\mathcal{A}(W_{11}^{12}\xi+b_1^2)),(W^{23}\mathcal{A}(W_{12}^{12}\xi+b_2^2))\right]\right\},\\
\tilde{W}_1^{I=2}\left(\xi;\mathbf{W},\mathbf{b},\mathcal{A}\right) &= W_{11}^{34}\left\{M\left[(W^{23}\mathcal{A}(W_{14}^{12}\xi+b_4^2)),(W^{23}\mathcal{A}(W_{15}^{12}\xi+b_5^2))\right]\right\},\\
\tilde{W}_2^{I=2}\left(\xi;\mathbf{W},\mathbf{b},\mathcal{A}\right) &= W_{22}^{34}\left\{M\left[(W^{23}\mathcal{A}(W_{11}^{12}\xi+b_1^2)),(W^{23}\mathcal{A}(W_{15}^{12}\xi+b_5^2))\right]\right\},\\
\tilde{W}_3^{I=2}\left(\xi;\mathbf{W},\mathbf{b},\mathcal{A}\right) &= W_{33}^{34}\left\{M\left[(W^{23}\mathcal{A}(W_{11}^{12}\xi+b_1^2)),(W^{23}\mathcal{A}(W_{14}^{12}\xi+b_4^2))\right]\right\}.
\end{aligned}
\tag{15}
$$

in which $M$ stands for the multiplication building block as defined in **Fig. 1**. Finally, the C-HiDeNN shape functions are given as

$$
\begin{aligned}
\tilde{N}_0(\xi;\mathbf{W},\mathbf{b},\boldsymbol{\mathcal{A}}) &= M\Big[(N_1(\xi;\mathbf{W},\mathbf{b},\boldsymbol{\mathcal{A}})),(\tilde{W}_0^{I=1}\big(\xi;\mathbf{W},\mathbf{b},\boldsymbol{\mathcal{A}})\big)\Big],\\
\tilde{N}_1(\xi;\mathbf{W},\mathbf{b},\boldsymbol{\mathcal{A}}) &= M\Big[(N_1(\xi;\mathbf{W},\mathbf{b},\boldsymbol{\mathcal{A}})),(\tilde{W}_1^{I=1}\big(\xi;\mathbf{W},\mathbf{b},\boldsymbol{\mathcal{A}})\big)\Big]\\
&\quad + M\Big[(N_2(\xi;\mathbf{W},\mathbf{b},\boldsymbol{\mathcal{A}})),(\tilde{W}_1^{I=2}\big(\xi;\mathbf{W},\mathbf{b},\boldsymbol{\mathcal{A}})\big)\Big],\\
\tilde{N}_2(\xi;\mathbf{W},\mathbf{b},\boldsymbol{\mathcal{A}}) &= M\Big[(N_1(\xi;\mathbf{W},\mathbf{b},\boldsymbol{\mathcal{A}})),(\tilde{W}_2^{I=1}\big(\xi;\mathbf{W},\mathbf{b},\boldsymbol{\mathcal{A}})\big)\Big]\\
&\quad + M\Big[(N_2(\xi;\mathbf{W},\mathbf{b},\boldsymbol{\mathcal{A}})),(\tilde{W}_2^{I=2}\big(\xi;\mathbf{W},\mathbf{b},\boldsymbol{\mathcal{A}})\big)\Big],\\
\tilde{N}_3(\xi;\mathbf{W},\mathbf{b},\boldsymbol{\mathcal{A}}) &= M\Big[(N_2(\xi;\mathbf{W},\mathbf{b},\boldsymbol{\mathcal{A}})),(\tilde{W}_3^{I=2}\big(\xi;\mathbf{W},\mathbf{b},\boldsymbol{\mathcal{A}})\big)\Big].
\end{aligned} \tag{16}
$$

Similar to the prior work on nonlinear HiDeNN-FEM implementation [12], three building blocks are introduced to evaluate the shape function derivative. Compared to the previous implementation, two key differences are highlighted here:

1. The partial derivative operator is carried over a convolution patch as opposed to just within the element.
2. In addition to optimizing the nodal position '$X$' (also known as $r$-adaptivity), C-HiDeNN achieves higher-order accuracy through the selection of polynomial order '$p$,' dilation parameter '$a$,' and patch size '$s$,'

In the 1D case, the shape function derivative matrix $\boldsymbol{\mathcal{B}}_{jI}^0 = \dfrac{\partial \tilde{N}_I}{\partial X_j}$ in Eq. (4) is simplified as $\boldsymbol{\mathcal{B}}_I^0 = \dfrac{d\tilde{N}_I}{dX}$ and is obtained through

$$
\begin{aligned}
\frac{du^h(\xi)}{dX} &= \sum_{I\in A_s^e}\frac{d\tilde{N}_I(\xi)}{dX}u_I = \sum_{I\in A^e}\frac{dN_I(\xi)}{dX}\sum_{J\in A_s^I}\tilde{W}_{a,p,J}(\xi)u_J + \sum_{I\in A^e}N_I(\xi)\sum_{J\in A_s^I}\frac{d\tilde{W}_{a,p,J}(\xi)}{dX}u_J\\
&= \sum_{I\in A_s^e}\boldsymbol{\mathcal{B}}_I^0\left(\xi\right)u_I
\end{aligned} \tag{17}
$$

Note that the steps to construct the neural network for the standard FE shape function derivative $\frac{dN_I(\xi)}{dX}$ in Eq.(17) follows the approach described in HiDeNN-FEM [12].

Derivatives of the convolution patch function $\widetilde{W}_{a,p,J}$ depends on the way it is constructed. In the case of 1D Lagrangian polynomial, the derivative can be analytically evaluated using the chain rule,

$$\frac{d\tilde{W}}{dX}=\frac{d\tilde{W}(\xi)}{d\xi}\left(\frac{dX}{d\xi}\right)^{-1}=\mathbf{J}^{-1}\frac{d\tilde{W}(\xi)}{d\xi} \tag{18}$$

in which Jacobian $\mathbf{J}=\dfrac{dX}{d\xi}$ is introduced. In the case of $p=2$, this is evaluated through

$$\mathbf{J}=\frac{\partial X}{\partial \xi}=\begin{bmatrix}\frac{d\tilde{W}_1}{d\xi} & \frac{d\tilde{W}_2}{d\xi} & \frac{d\tilde{W}_3}{d\xi}\end{bmatrix}\begin{Bmatrix}X_1\\ X_2\\ X_3\end{Bmatrix}=\mathbf{D}_N\mathbf{X}. \tag{19}$$

where $\mathbf{X}=\{X_1 \quad X_2 \quad X_3\}^T$ provides the physical coordinates of the nodes.

**Fig. 7** (a) provides the neural network representations of the derivatives of the 1D convolution patch functions. The corresponding expressions are then given as

$$\begin{aligned}
\frac{d\tilde{W}_1}{dX}&=M\left[W_{11}^{67}\mathbf{J}^{-1},W_{11}^{37}\boldsymbol{\mathcal{A}}(W_{11}^{23}\boldsymbol{\mathcal{A}}(W_{11}^{12}\xi+b_1^2))\right],\\
\frac{d\tilde{W}_2}{dX}&=M\left[W_{12}^{67}\mathbf{J}^{-1},W_{22}^{37}\boldsymbol{\mathcal{A}}(W_{22}^{23}\boldsymbol{\mathcal{A}}(W_{12}^{12}\xi+b_2^2))\right],\\
\frac{d\tilde{W}_3}{dX}&=M\left[W_{13}^{67}\mathbf{J}^{-1},W_{33}^{37}\boldsymbol{\mathcal{A}}(W_{33}^{23}\boldsymbol{\mathcal{A}}(W_{13}^{12}\xi+b_3^2))\right].
\end{aligned} \tag{20}$$

with

$$\mathbf{J}=W_{11}^{56}\boldsymbol{\mathcal{A}}\begin{Bmatrix}W_{11}^{45}\boldsymbol{\mathcal{A}}(W_{11}^{34}\boldsymbol{\mathcal{A}}(W_{11}^{23}\boldsymbol{\mathcal{A}}(W_{11}^{12}\xi+b_1^2)))+\\ W_{21}^{45}\boldsymbol{\mathcal{A}}(W_{22}^{34}\boldsymbol{\mathcal{A}}(W_{22}^{23}\boldsymbol{\mathcal{A}}(W_{12}^{12}\xi+b_2^2)))+\\ W_{31}^{45}\boldsymbol{\mathcal{A}}(W_{33}^{34}\boldsymbol{\mathcal{A}}(W_{33}^{23}\boldsymbol{\mathcal{A}}(W_{13}^{12}\xi+b_3^2)))\end{Bmatrix} \tag{21}$$

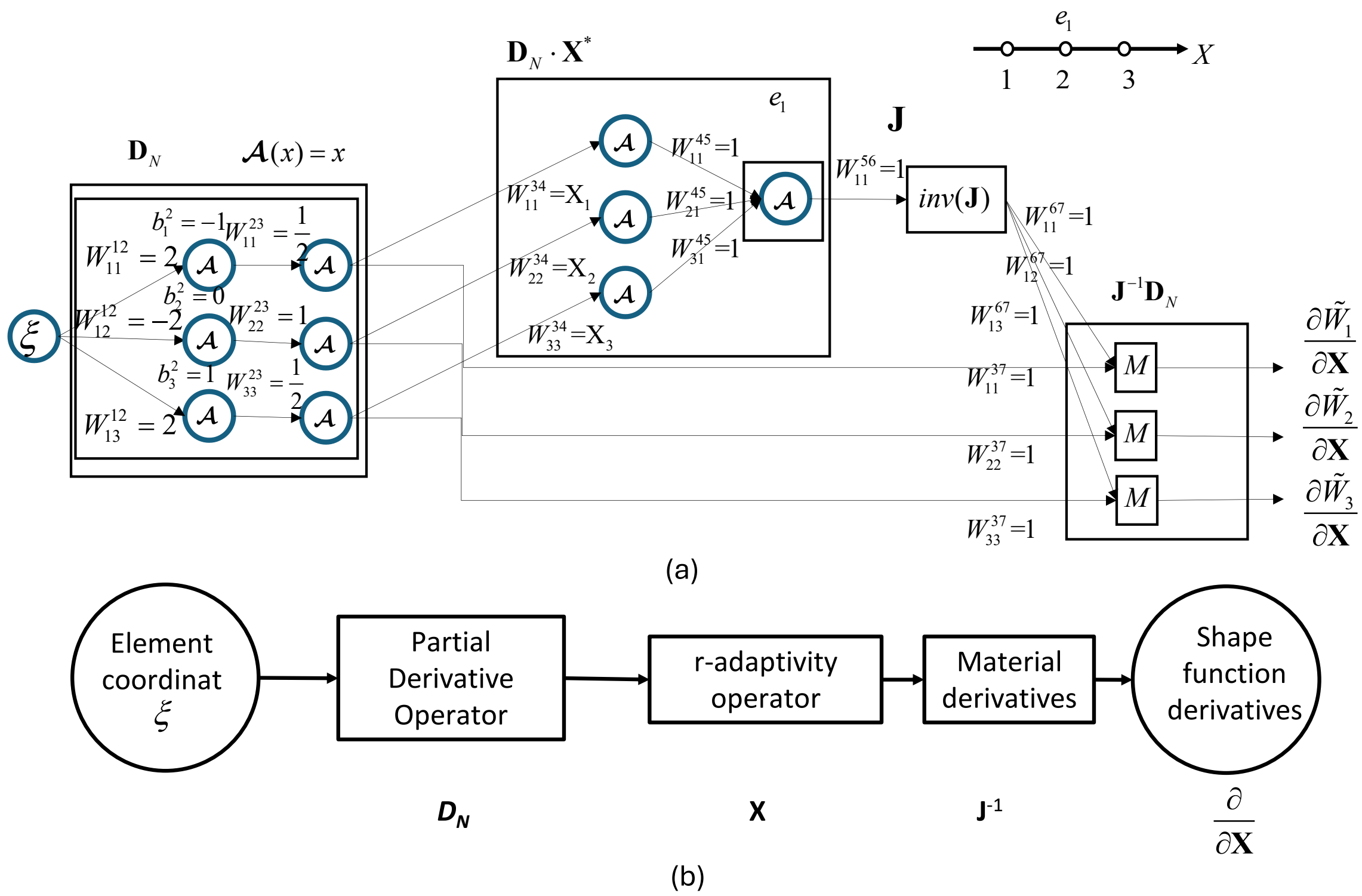

**Fig. 7. (a) Neural network representation for evaluating the derivative of the convolution kernel in 1D. (b) A general flowchart for evaluating the shape function derivatives.**

The steps outlined above can be extended to higher dimensional cases according to **Fig.7b**. The general methodology to construct NN representation for the C-HiDeNN shape function is illustrated in **Fig. 8**. In this implementation, the element and physical coordinates serve as inputs for generating the finite element (FE) shape functions and convolution patch functions, respectively. The process of generating FE shape functions follows the HiDeNN-FE approach and has been extensively discussed in ref [11, 12]. The construction of the NN for the convolution patch functions, on the other hand, depends on the specific choice of the meshfree approximation employed.

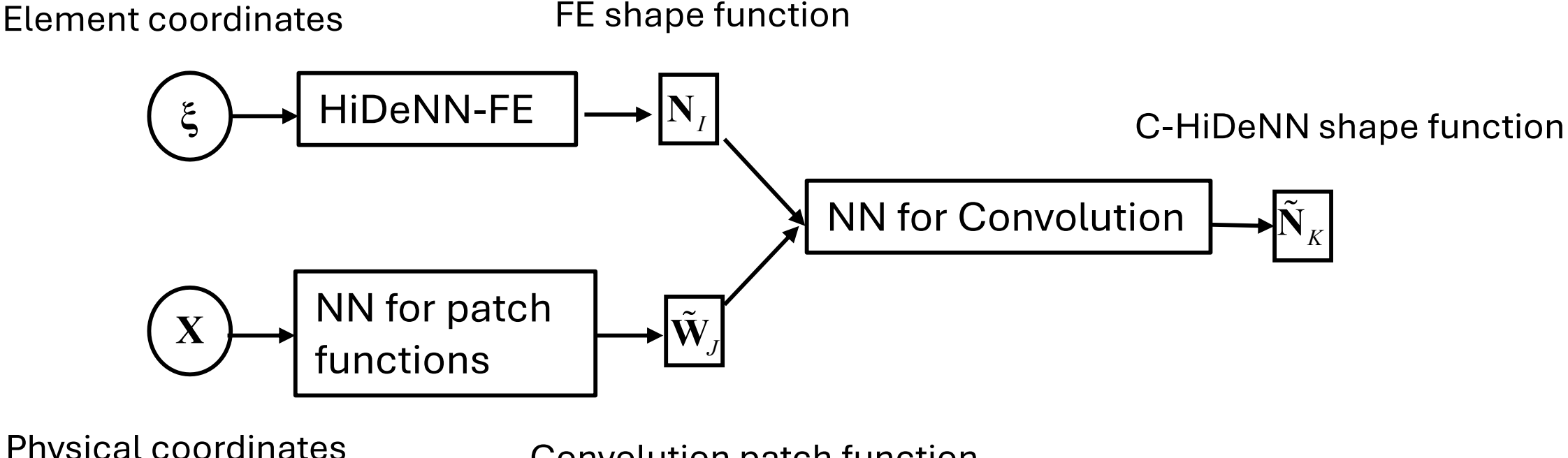

**Fig. 8. General approach to construct C-HiDeNN shape functions.**

Below we will demonstrate this using the 2D case illustrated in **Fig. 4** with the patch parameters $s$=1, $p$=1 and linear basis ( $m=3$ ). We consider the patch function for node 1. In this case the nodal patch is given as $A_s^{i=1}=\{1,2,3,4,11,12,13,14,15\}$ . Following Eq.(11), it can be given that

$$\mathbf{G}=\begin{bmatrix} B_1(\mathbf{X}_1) & B_2(\mathbf{X}_1) & \dots & B_9(\mathbf{X}_1) & 1 & X_1 & Y_1 \\ B_1(\mathbf{X}_2) & B_2(\mathbf{X}_2) & \dots & B_9(\mathbf{X}_2) & 1 & X_2 & Y_2 \\ \vdots & \vdots & & \vdots & \vdots & \vdots & \vdots \\ B_1(\mathbf{X}_9) & B_2(\mathbf{X}_9) & \dots & B_9(\mathbf{X}_9) & 1 & X_9 & Y_9 \\ & & & & & & \\ 1 & 1 & \dots & 1 & 0 & 0 & 0 \\ X_1 & X_2 & \dots & X_9 & 0 & 0 & 0 \\ Y_1 & Y_2 & \dots & Y_9 & 0 & 0 & 0 \end{bmatrix} \tag{22}$$

in which $B_I(\mathbf{X}_J)$ is the radial base function that centers at node $I$ and evaluated at node $J$ . Here $I,J=1,2,3,...,9$ is introduced as local nodal indices that map to the nodal patch of $A_s^{i=1}=\{1,2,3,4,11,12,13,14,15\}$ . Based on Eq.(22), element of $\mathbf{G}$ can be formed by constructing the NN that uses the nodal coordinates of the patch nodes as the input with the appropriate activation functions. For instance, $G_{29}=B_2(\mathbf{X}_9)$ is obtained by using $\mathbf{X}_9$ as the input and activation function of $B_2(\mathbf{X})=B(|\mathbf{X}-\mathbf{X}_2|)$ in which $r=|\mathbf{X}-\mathbf{X}_2|$ gives the radial distance. Once $\mathbf{G}$ is formed, its inverse can be obtained by $\mathbf{G}^{-1} = \frac{1}{\det(\mathbf{G})}\mathbf{G}_c^T$ in which $\mathbf{G}_c$ is a co-factor matrix of $\mathbf{G}$ and superscript "$T$" denotes transpose. As can be seen, this inverse operation involves algebraic operations of rational functions. A hierarchical NN can be constructed to realize the inverse operation using the fundamental building blocks of linear function, multiplication and inversion [11]. Finally, the patch functions $\tilde{\mathbf{W}}(\mathbf{X})$ and its derivatives are obtained from Eq. (10) and Eq. (12) in which the matrix multiplication can also be obtained constructing the NN.

## 3 Solution Scheme

Solving the discretized momentum equation from Eq.(3) is equivalent to the energy minimization statement of

$$\arg\min_{\mathbf{d}} \Pi(\mathbf{d}) \tag{23}$$

in which $\Pi = W_{kin} + W_{int} - W_{ext}$ and $W_{kin}$ , $W_{int}$ and $W_{ext}$ represent the kinetic, internal and external energy of the system, respectively. $\mathbf{d}$ is the nodal displacement vector and the displacement field is approximated by

$$\mathbf{u}^h = \tilde{\mathbf{N}}(\xi; \mathbf{X}, \mathbf{W}, \mathbf{b})\mathbf{d} \tag{24}$$

with $\tilde{\mathbf{N}}$ being the C-HiDeNN shape function matrix, $\xi$ and $\mathbf{X}$ are the element and physical coordinates, $\mathbf{W}$ and $\mathbf{b}$ are the associated weights and biases for the NN representation.

In dynamic simulation, it is a common practice to divide the total time domain into multiple intervals and solve the problem using an incremental solution scheme. This approach is typically coupled with a time discretization method, e.g., the central difference (CD) scheme. Assuming the solution $\mathbf{d}_n$ to Eq.(23) is known at time $t_n$, the goal is to find the solution at $t_{n+1} = t_n + \Delta t$. In CD, a half-time step is introduced as $t_{n+1/2} = t_n + \Delta t / 2$ where the velocity is given by $\mathbf{v}^{n+1/2}$ and approximated as

$$\mathbf{v}^{n+1/2} = \mathbf{v}^n + \left(t_{n+1/2} - t_n\right)\mathbf{a}^n \tag{25}$$

and

$$\mathbf{v}^{n+1} = \mathbf{v}^{n+1/2} + \left(t_{n+1} - t_{n+1/2}\right)\mathbf{a}^{n+1} \tag{26}$$

Here $\mathbf{a}^n$ and $\mathbf{a}^{n+1}$ are the acceleration at $t_n$ and $t_{n+1}$, respectively, and they are evaluated according to $\mathbf{a}^n = \mathbf{M}^{-1}\left(\mathbf{f}_{ext}^n - \mathbf{f}_{int}^n\right)$ and $\mathbf{a}^{n+1} = \mathbf{M}^{-1}\left(\mathbf{f}_{ext}^{n+1} - \mathbf{f}_{int}^{n+1}\right)$.

The incremental energy terms are then expressed as

$$\Delta W_{int} = W_{int}^{n+1} - W_{int}^n = \frac{1}{2}\Delta\mathbf{d}^T\left(\mathbf{f}_{int}^n + \mathbf{f}_{int}^{n+1}\right) \tag{27}$$

$$\Delta W_{ext} = W_{ext}^{n+1} - W_{ext}^n = \frac{1}{2}\Delta\mathbf{d}^T\left(\mathbf{f}_{ext}^n + \mathbf{f}_{ext}^{n+1}\right) \tag{28}$$

$$\Delta W_{kin} = W_{kin}^{n+1} - W_{kin}^n = \frac{1}{2}\left(\mathbf{v}^{n+1}\right)^T \mathbf{M}\left(\mathbf{v}^{n+1}\right)^T - \frac{1}{2}\left(\mathbf{v}^n\right)^T \mathbf{M}\left(\mathbf{v}^n\right)^T \tag{29}$$

Based on Eq. (27) to (29), the energy minimization statement of Eq.(23) is revised to be

$$\arg\min_{\Delta\mathbf{d}}\left(\Pi\left(\mathbf{d}^n\right) + \Delta\Pi\left(\Delta\mathbf{d}\right)\right) \tag{30}$$

with $\Delta\Pi = \Delta W_{kin} + \Delta W_{int} - \Delta W_{ext}$. For each time increment, Eq. (30) is solved to update the displacement solution based on $\Delta\Pi$ and Eqs. (27) to (29).

## 4 Numerical Examples

In this section, we provide several numerical examples to demonstrate the performance of C-HiDeNN. An in-house code has been developed for this purpose. The discretized nonlinear momentum equation (Eq.(3)) is solved using solution algorithm described in section 3. All computations were performed on a workstation with Intel(R) Core(TM) i7-9750H CPU @ 2.60GHz and 32 GB memory.

### 4.1 A plate with a notch subjected to dynamic load and unload condition

We consider a plane strain problem involving a notched rectangular plate subjected to combined horizontal and vertical displacements applied to the top surface, as illustrated in **Fig. 9**. The plate measures 1 m by 0.3 m, with a notch radius of R=0.1m. The displacement history is presented in **Fig. 9** (d). and the plate is fixed at the bottom. **Fig. 9** (b) shows the finite element mesh, consisting of 450 4-node quadrilateral elements, which is also used for generating the C-HiDeNN shape functions and performing the analysis. For comparison, **Fig. 9** (c) presents the mesh used for the reference solution solved using ABAQUS. This model includes 533,422 CPE4R elements (four-node plane strain elements with reduced integration), resulting in 1,070,628 degrees of freedom. The plate is modeled as Neo-Hookean material with the strain energy density function *w* given as

$$w=C_{10}(\overline{I}_1-3)+\frac{1}{D_1}(J\text{-}1)^2, \tag{31}$$

where, $\overline{I}_1=J^{-2/3}I_1$ , $J$ = det(**F**) and $I_1$ is the first invariant of the right Cauchy-Green deformation tensor. The material parameters are given as $C_{10}$ = 115.385 kPa, $D_1$ = 4×10$^{-6}$ Pa$^{-1}$.

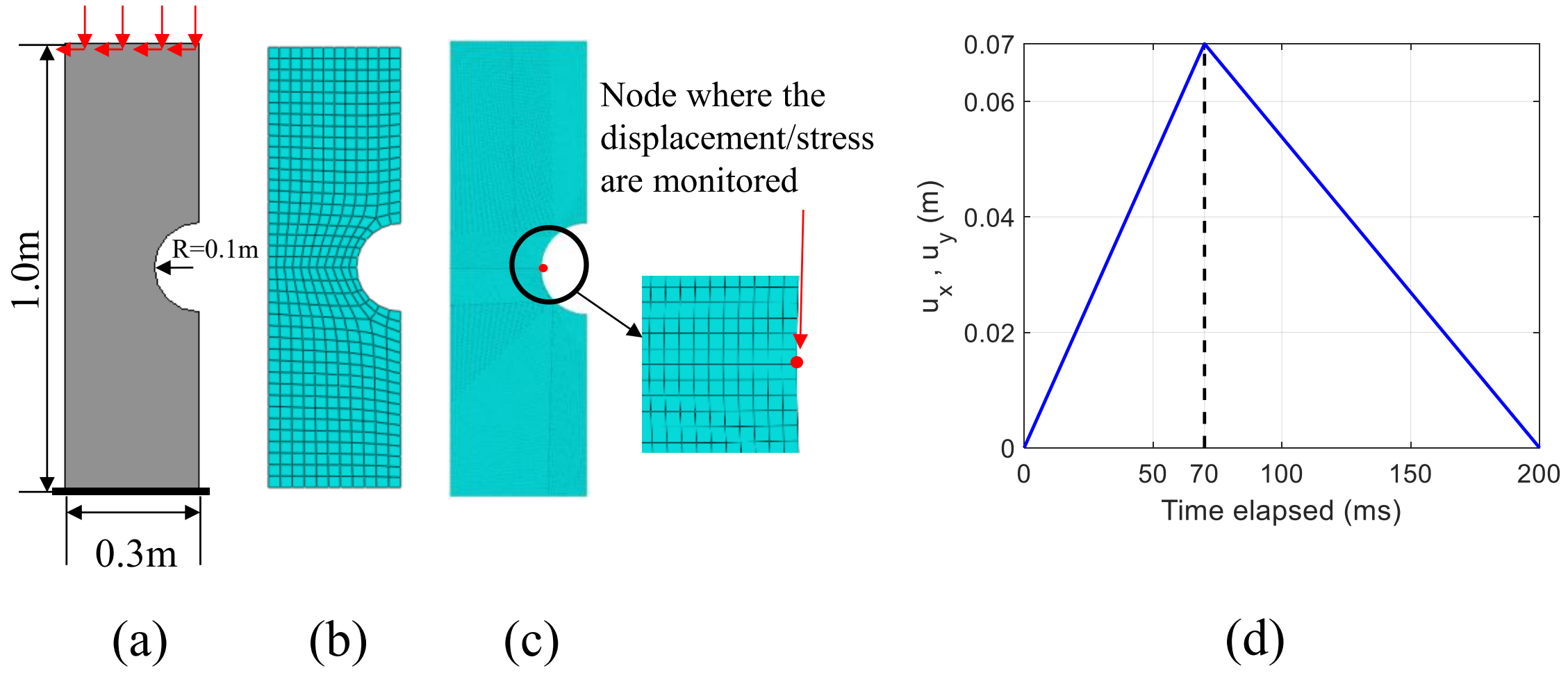

**Fig. 9. A plate with a notch subjected to load and unload condition. (a) Dimension of the plate with notch, (b) The finite element mesh (also used for C-HiDeNN), (c) Reference mesh (CPE4R), (d) History of displacements applied on the top surface.**

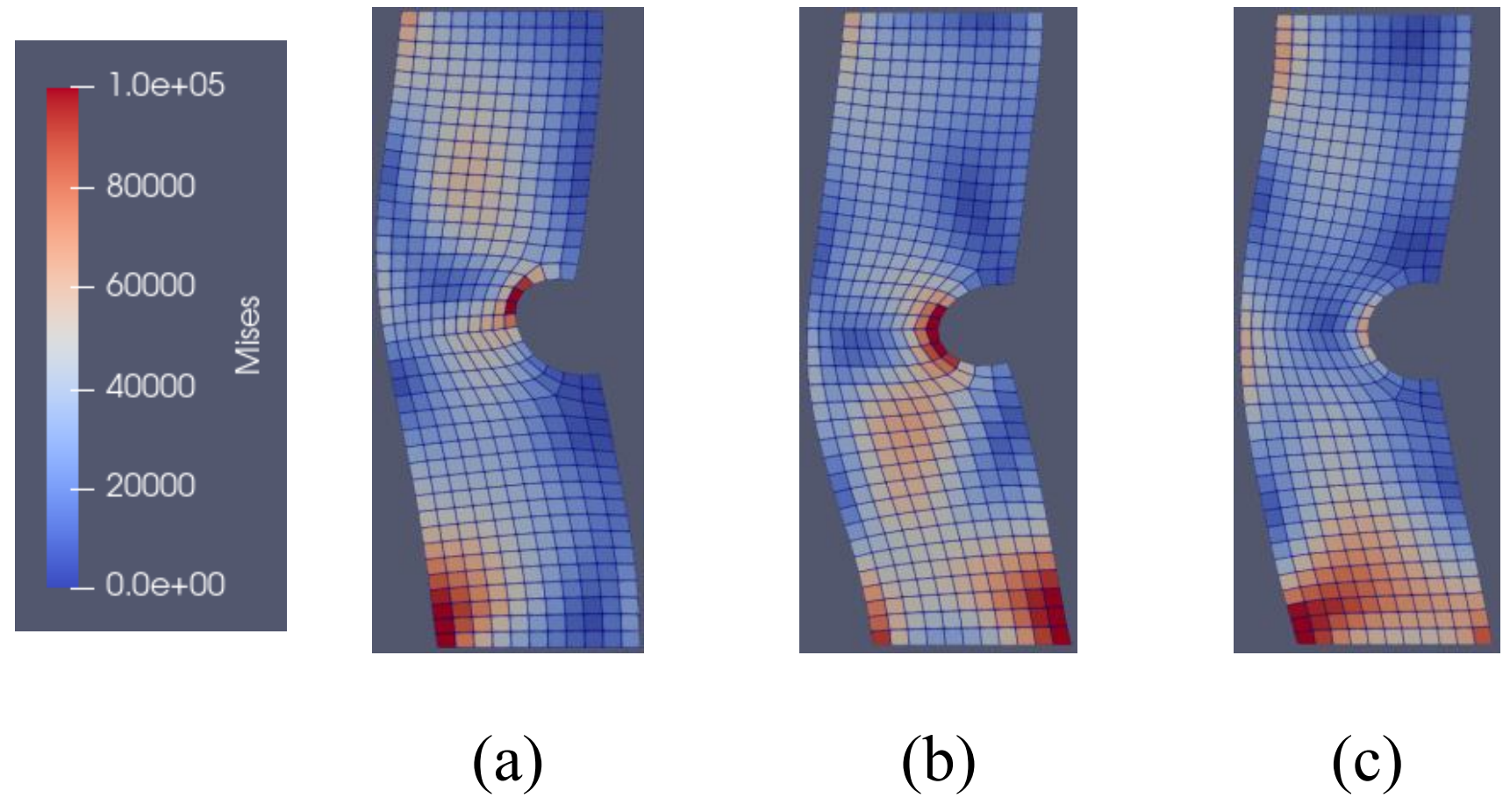

**Fig. 10. von Mises stress (Pa) distribution at different time steps: (a) 100ms, (b) 120ms, (c) 140ms.**

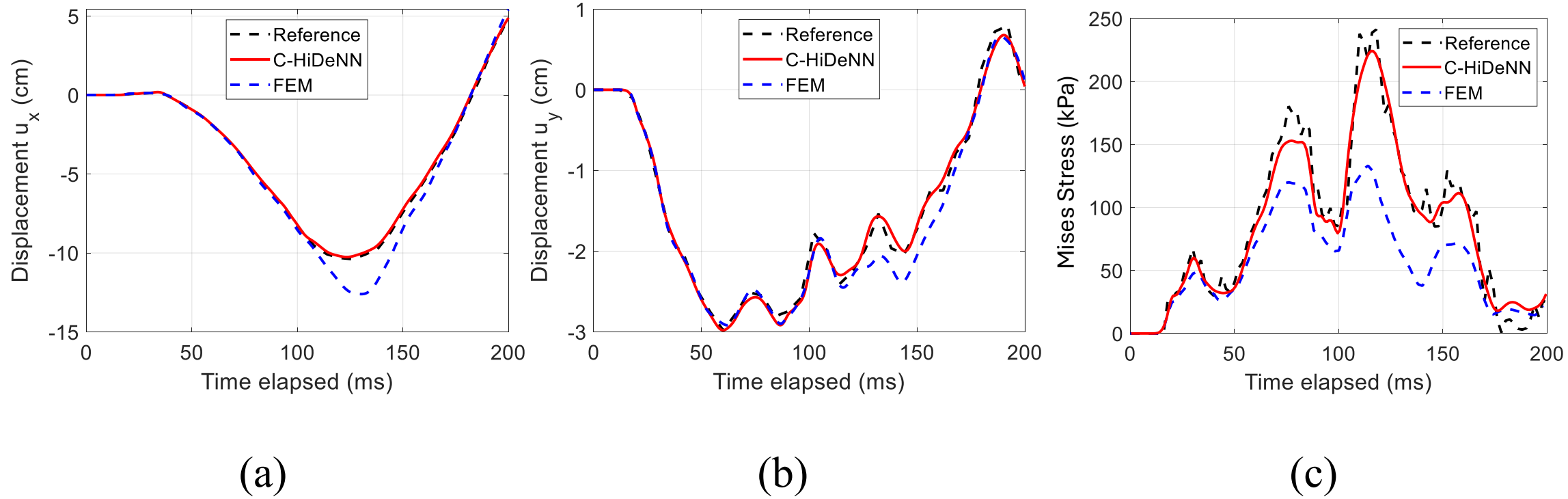

**Fig. 11. Comparison of the solution at the notch tip using reference mesh with** 533,422 **CPE4R elements, C-HiDeNN, and standard FEM: (a) x-axis displacement** $u_x$ **(cm), (b) y-axis displacement** $u_y$ **(cm), (c) von Mises stress (kPa).**

| Analysis/field result | FEM (CPE4R) | C-HiDeNN |
|---|---|---|
| $u_x$ | 23.4% | 2.4% |
| $u_y$ | 37.5% | 2.2% |
| $\sigma_{mises}$ | 50.7% | 4.0% |

**Table 1: Error in the displacements and effective stress at 130ms.**

The stress contour figure at different times is shown in **Fig. 10**. The results at the monitored notch point are shown in **Fig. 11** and the error analysis result is shown in **Table 1**. Results show that C-HiDeNN method can capture displacement history (**Fig. 11** (a) $u_x$ and (b) $u_y$) and (c) the stress concentration, much better than finite element method with the same mesh. Error analysis at 130ms in **Table 1** shows that C-HiDeNN acquires the displacement and stress solution with higher accuracy (less than 5%) than FEM (over 23%).

In addition to the full-domain approach described above, C-HiDeNN can also be applied selectively to regions requiring high resolution, without modifying the finite element (FE) mesh. We refer to the method known as *s*-adaptivity. In this example, *s*-adaptivity is implemented by introducing the C-HiDeNN approximation in the area surrounding the notch, as illustrated in **Fig. 12** (a), while 4-node quadrilateral elements are used for the remainder of the domain. The prediction results, shown in **Fig. 12**, are compared with both the reference solution and the FE solution using the same mesh. As shown, the *s*-adaptivity approach maintains high accuracy while significantly reducing computational time, from 634 seconds with full-domain C-HiDeNN to just 56 seconds with the *s*-adaptivity implementation.

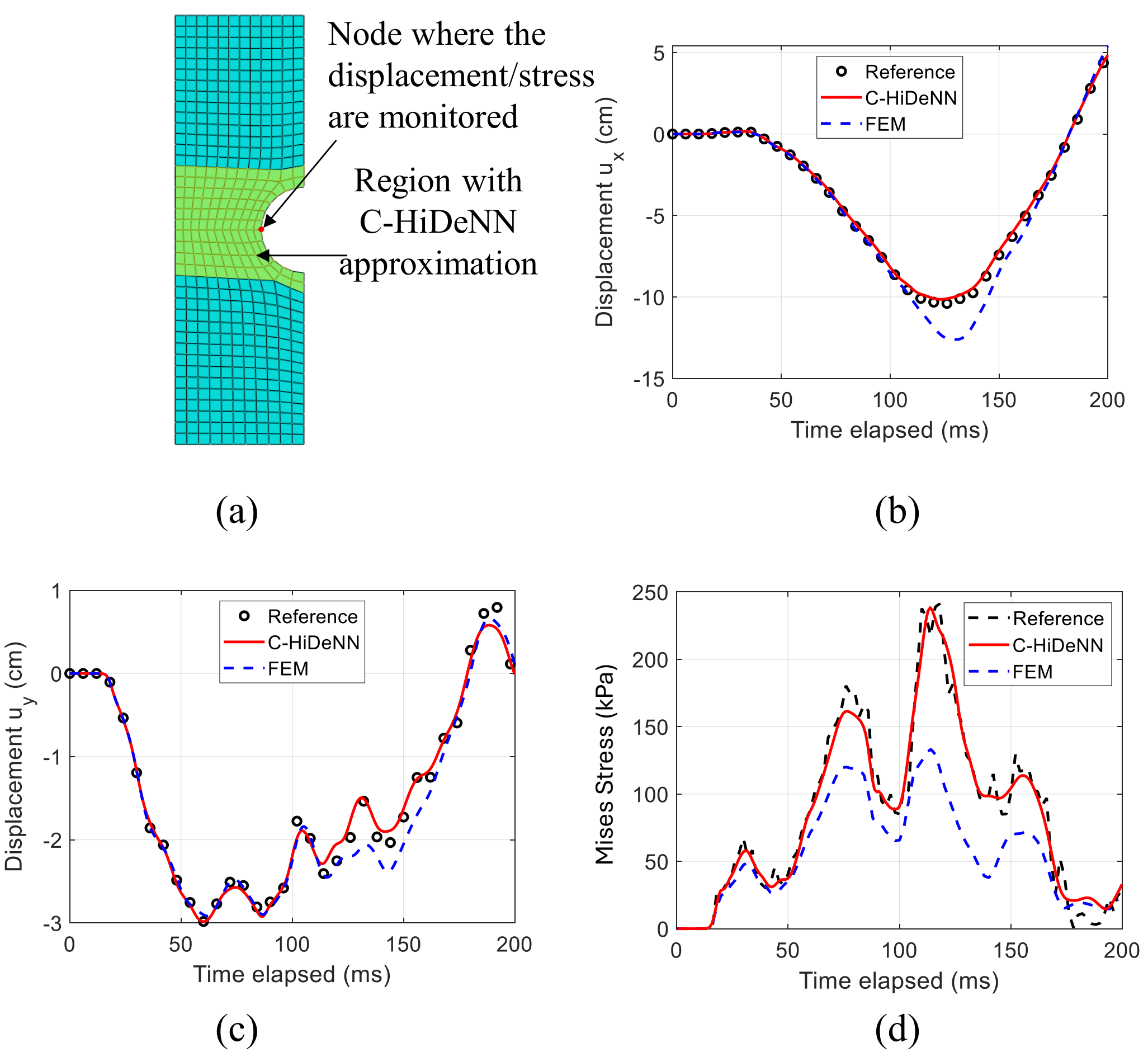

**Fig. 12. (a) Hybrid FEM/C-HiDeNN discretization. Comparison of the results at the monitored node near the notch tip among reference solution (CPE4R), C-HiDeNN, and standard FEM: (b) x-axis displacement $u_x$ (cm), (c) y-axis displacement $u_y$ (cm), (d) von Mises stress (kPa).**

### 4.2 A plate with a notch subjected to cyclic load

In this case, a plate with a notch [30] subjected to cyclic loading is considered. The plate is modeled using the same Neo-Hookean material as in Case 1. The sample measures 0.5 m by 1 m, with a notch radius of R=0.05 m, as shown in **Fig. 13** (a). **Fig. 13** (b) shows the finite element mesh consisting of 1,200 four-node quadrilateral elements, which is also used for subsequent C-HiDeNN analysis. The reference solution is obtained using a highly refined mesh with 510,101 CPE4R elements and 1,020,202 degrees of freedom. This reference mesh provides the converged solution against which results from the coarse mesh, computed using both standard FEM and the C-HiDeNN method, are compared. For the boundary conditions, the bottom surface of the model is fixed, and a pressure load is applied to the top surface according to the load history shown in **Fig. 13** (c). The displacement and stress results at the center of the notch are presented in **Fig. 14**. As shown in **Fig. 14** (a) for $u_x$ and (b) $u_y$, , and **Fig. 14** (c) for stress, the C-HiDeNN method captures the displacement history and stress distribution significantly better than the finite element method using the same mesh.

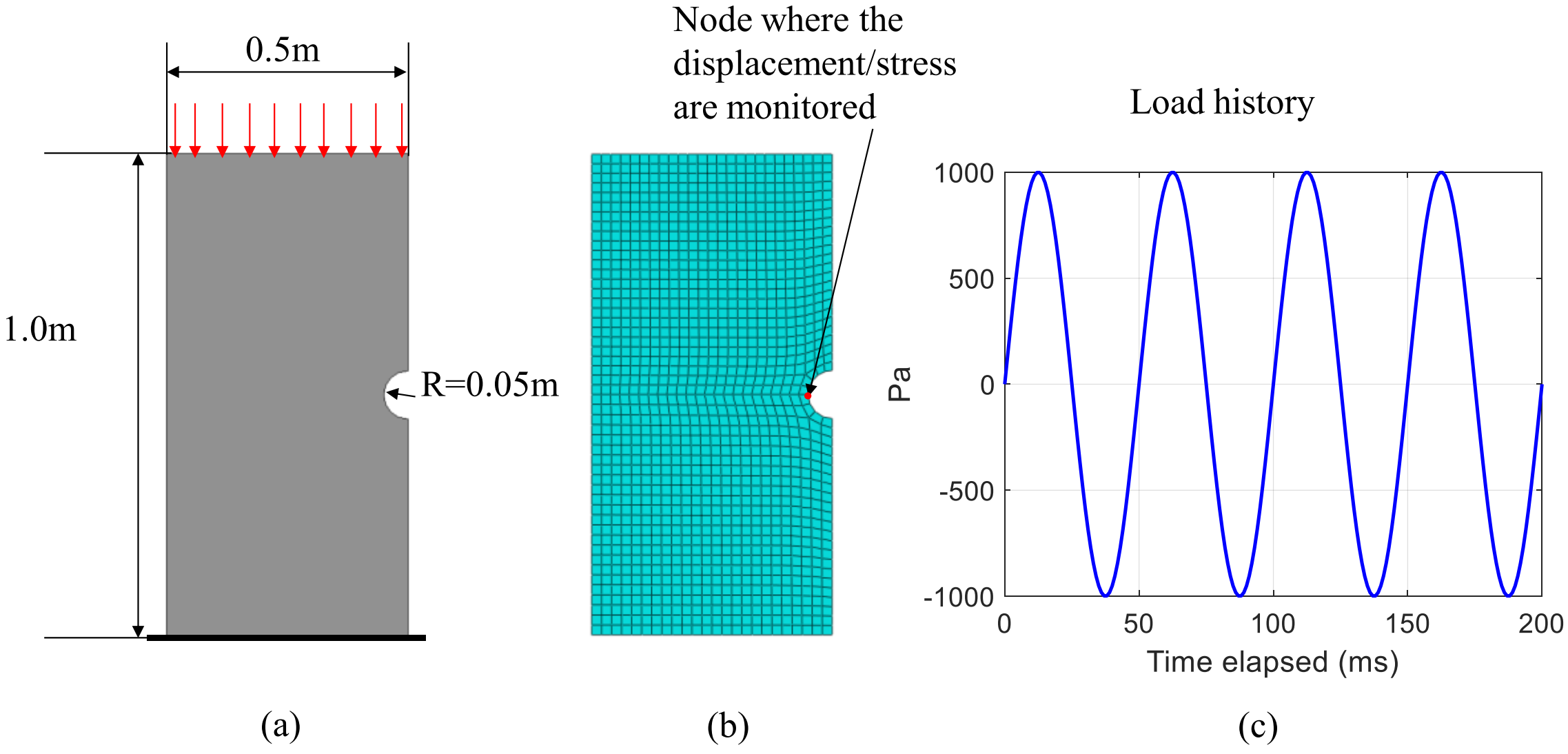

**Fig. 13. A plate with a notch subjected to cyclic load (a) Dimensions of the plate with notch. (b) The finite element mesh (also for C-HiDeNN). (c) Time history of the applied cyclic load.**

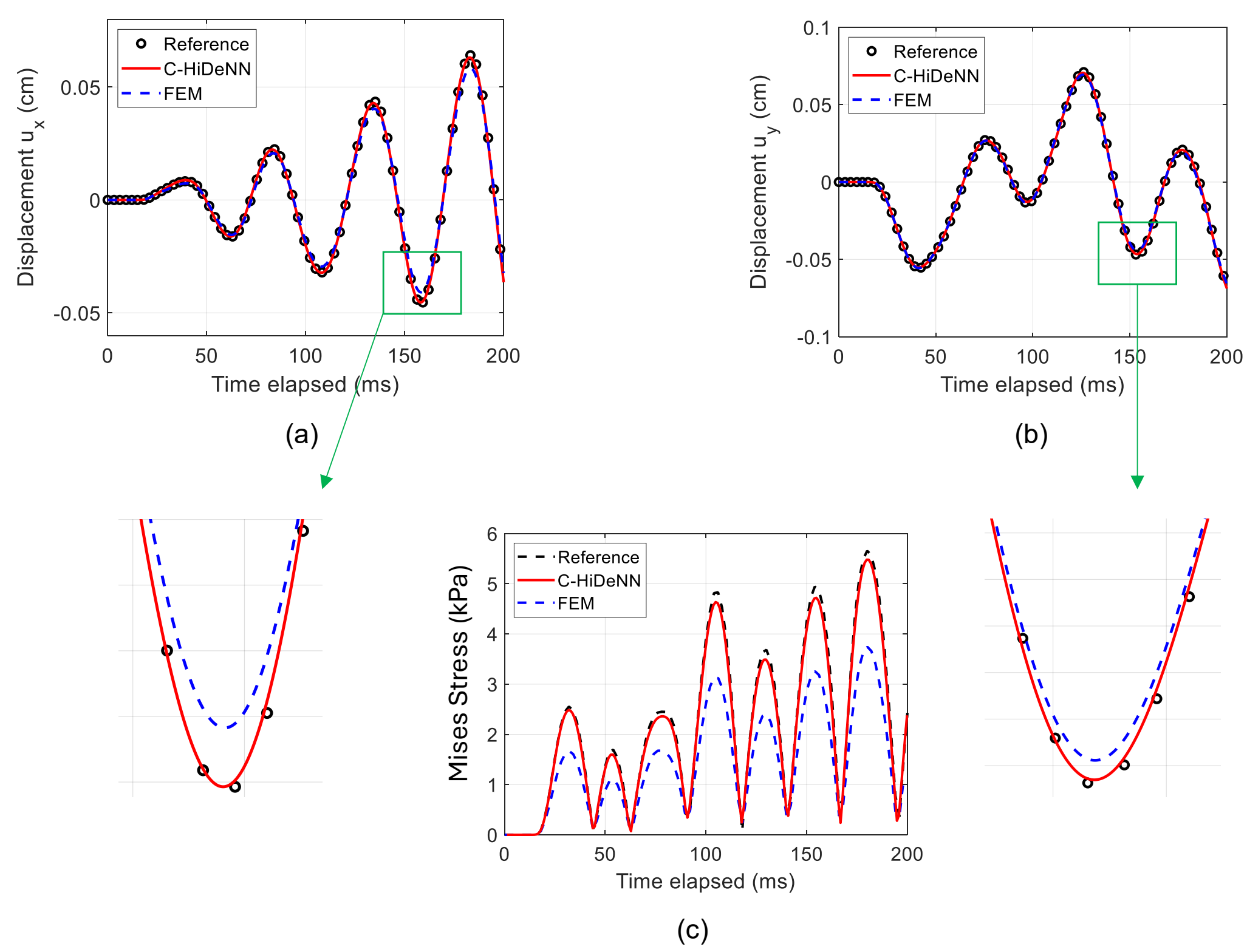

**Fig. 14. Comparison on the results near the notch tip among three methods: reference solution (CPE4R), C-HiDeNN and standard FEM. (a) Displacement in the x-direction ($u_x$, cm), (b) Displacement in the y-direction ($u_y$, cm), (c) von Mises stress (kPa).**

Similar to the last case, we also implemented *s*-adaptivity in which only the center region around the notch is approximated with C-HiDeNN, as shown in **Fig. 15** (a). The prediction results (displacement and stress) are compared with the full field implementation and reference solution and shown in **Fig. 15**. The *s*-adaptive approach was found to maintain high accuracy. In terms of computational efficiency, the total runtime was reduced from 4.3 hours for the reference solution and 742 seconds for the full-field C-HiDeNN case to just 82 seconds with the C-HiDeNN *s*-adaptivity implementation.

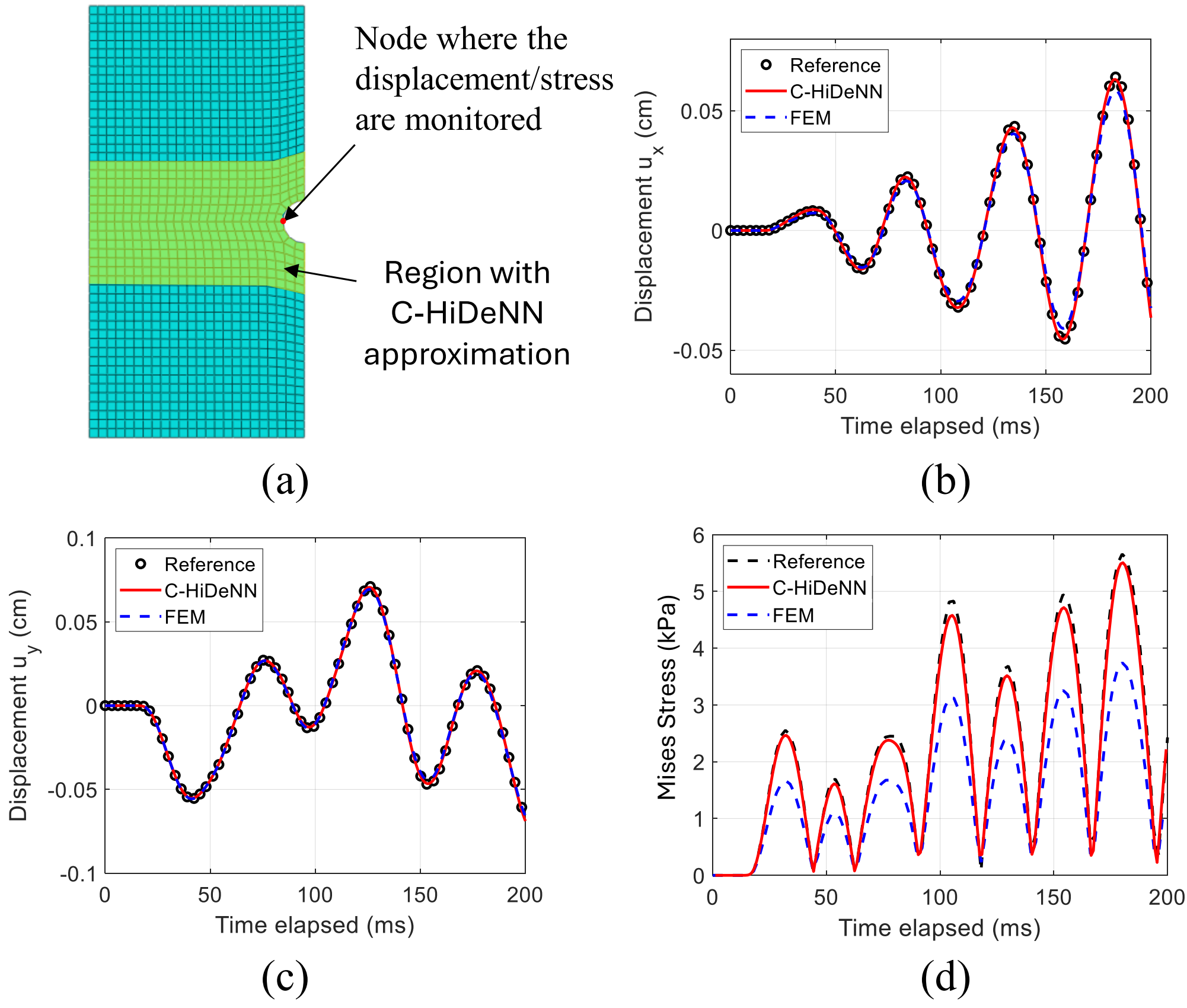

**Fig. 15. (a) Hybrid FEM/C-HiDeNN discretization. Comparison on the results at the monitored node near the notch tip among reference solution (CPE4R), C-HiDeNN, and standard FEM: (b) Displacement in the x-direction ($u_x$, cm), (c) Displacement in the y-direction ($u_y$, cm), (d) von Mises stress (kPa).**

### 4.3 A 2D model of mouse brain under impact

Mouse brain is frequently used as an experimental model to investigate the mechanics and effects of traumatic brain injuries [31]. In this study, we analyze a 2D cross-sectional geometry of a mouse brain model that was established based on 2D mouse brain image [31], treating it as a plane strain problem. The geometry is depicted in **Fig. 16** (a), where the bottom two regions are fixed, and an impact load is applied to the top surface over a duration of 5 milliseconds, as shown in **Fig. 16** (b). The mouse brain is modeled as neo-Hookean material with two sets of parameters introduced (Materials 1 and 2 as shown in **Fig. 16**). Material 1 has an initial shear modulus 1730.8Pa, and an initial bulk modulus 3750Pa. Material 2 has an initial shear modulus 4230.8Pa, and an initial bulk modulus 916Pa. The reference solution is obtained using a fine mesh consisting of 516,558 CPE4R elements and 1,035,956 degrees of freedom. For comparison, the same problem is also analyzed using FEM and C-HiDeNN, both employing a coarse mesh of 1,162 quadrilateral elements with 2,474 degrees of freedom. C-HiDeNN requires 334 seconds of computational time, compared to 4.7 hours for the reference solution.

The von Mises stress distribution (in Pascals) at 2 ms, 3 ms, and 4 ms is shown in **Fig. 17**. To evaluate the dynamic response, a node located at the center of the mouse brain is selected to monitor displacement and stress over time. The results are presented in **Fig. 18**. As shown in **Fig. 18** (a) for $u_x$ and (b) for $u_y$, both FEM and C-HiDeNN capture the displacement history with reasonable accuracy. However, C-HiDeNN demonstrates significantly better accuracy in capturing the von Mises stress within the brain region, as shown in **Fig. 18** (c).

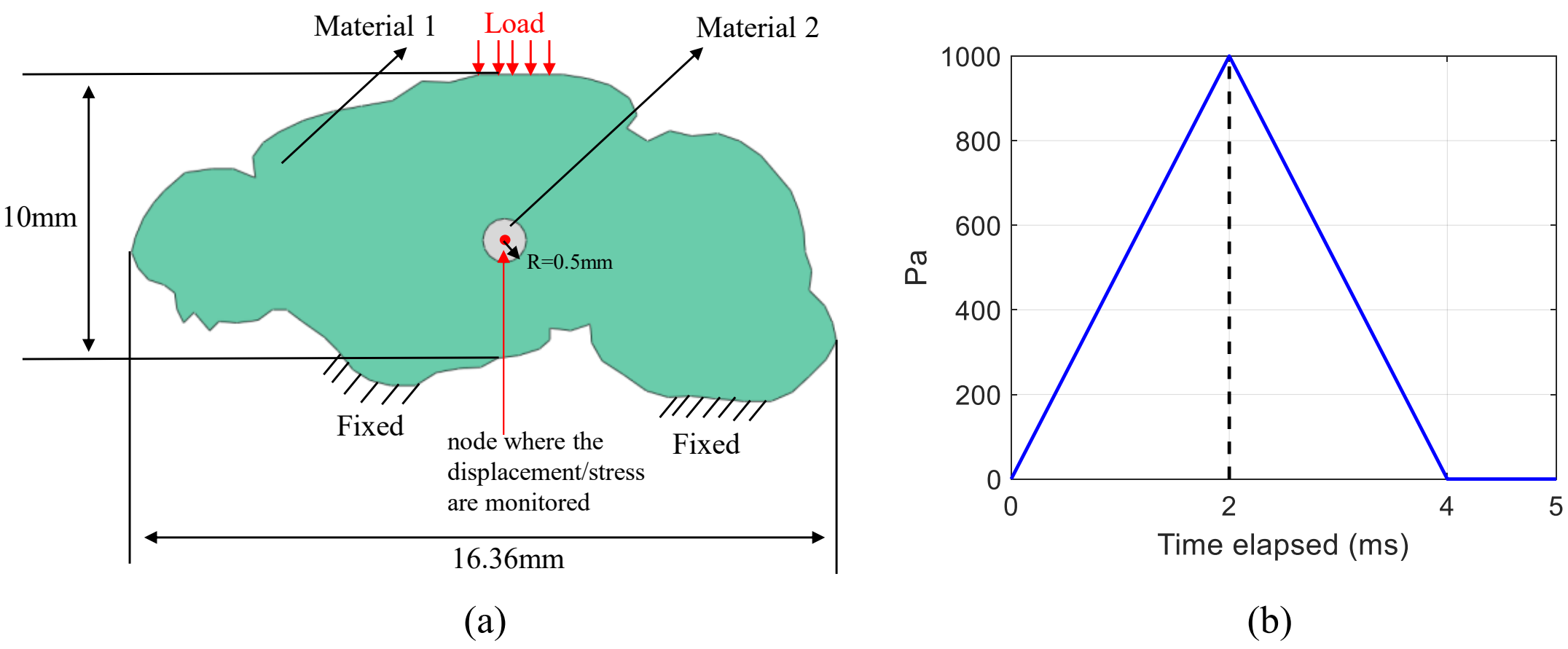

**Fig. 16. Simulation of mouse brain under impact: (a) geometry, dimensions and boundary conditions, (b) History of the impact load.**

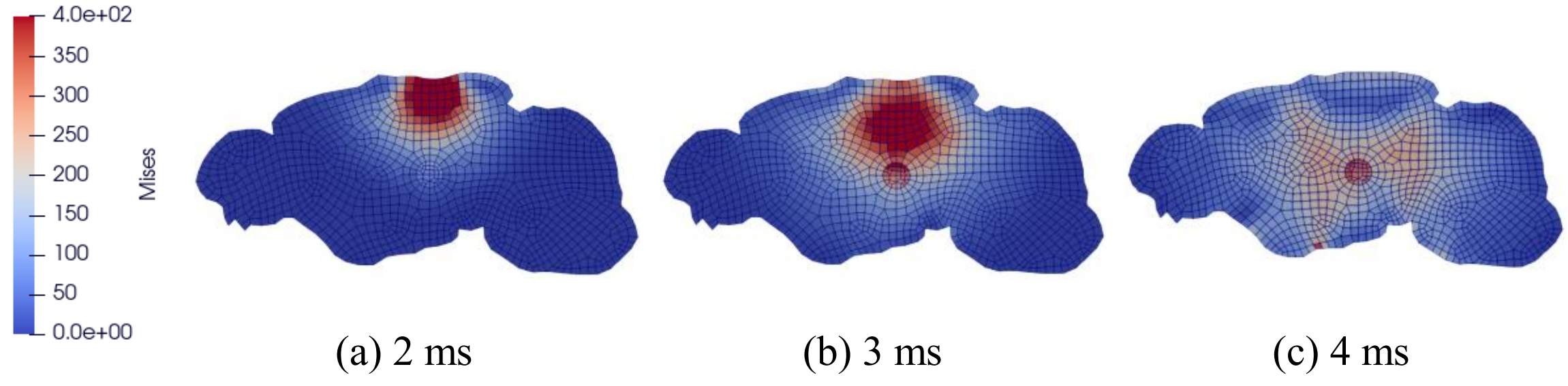

**Fig. 17. von Mises stress (Pa) of the mouse brain at different time steps: (a) 2ms, (b) 3ms, (c) 4ms.**

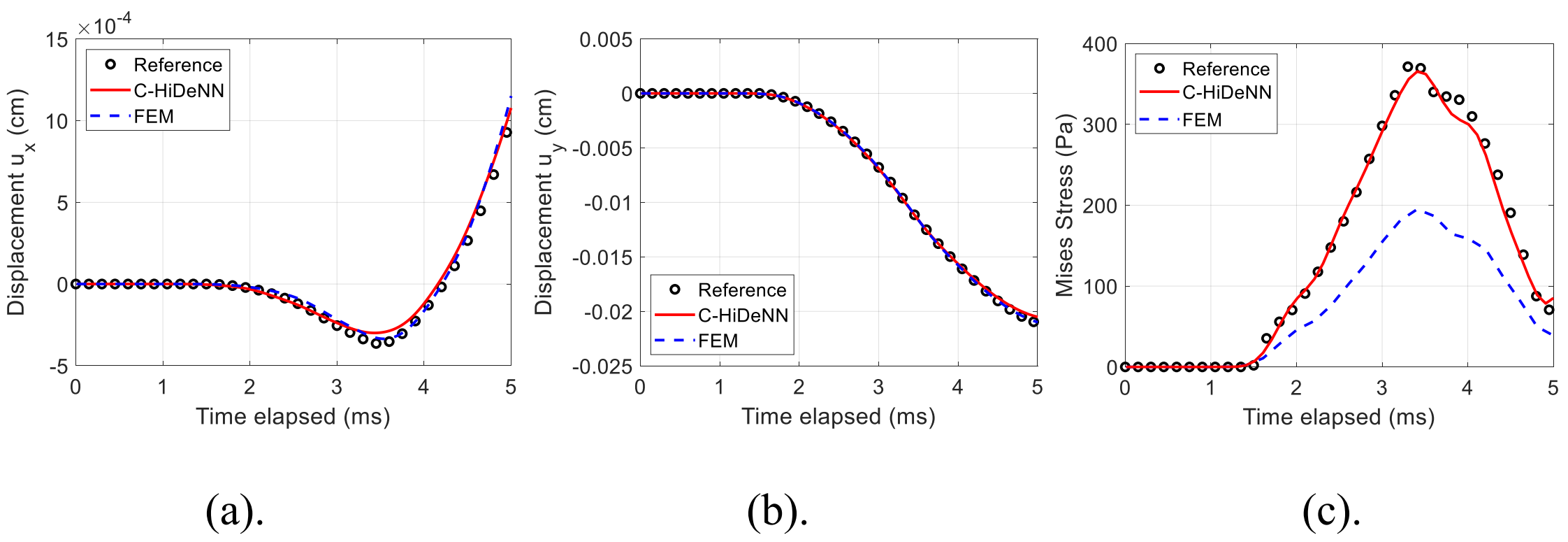

**Fig. 18. Comparison of results at the monitored point for the mouse brain impact problem among the reference solution, C-HiDeNN, and standard FEM: (a) Displacement in the x-direction ($u_x$, cm), (b) Displacement in the y-direction ($u_y$, cm), (c) von Mises stress (Pa).**

### 4.4 A plate with a long notch under tensile load condition

We consider a 2D plane strain simulation of a plate with a long notch subjected to the pressure load as shown in **Fig. 19**. The sample dimensions are 1.0 m by 1.0 m, with a notch measuring 25 cm by 2 cm located at the center of the right-hand edge. The material model is the same as in the previous case. The bottom surface is fixed, and a displacement load is applied to the top surface. The loading history is shown in **Fig. 20**, with a 1 cm displacement applied over 40 ms.

The reference solution is obtained using a fine mesh of 5,127,370 CPE4R elements, resulting in 10,264,992 degrees of freedom. The same problem is analyzed using both FEM and C-HiDeNN with a coarse mesh of 4,424 quadrilateral elements and 9,154 degrees of freedom. C-HiDeNN requires 847 seconds of computational time, compared to 3.4 hours for the reference solution. Stress distribution contours at 20 ms and 30 ms are shown in **Fig. 21**.

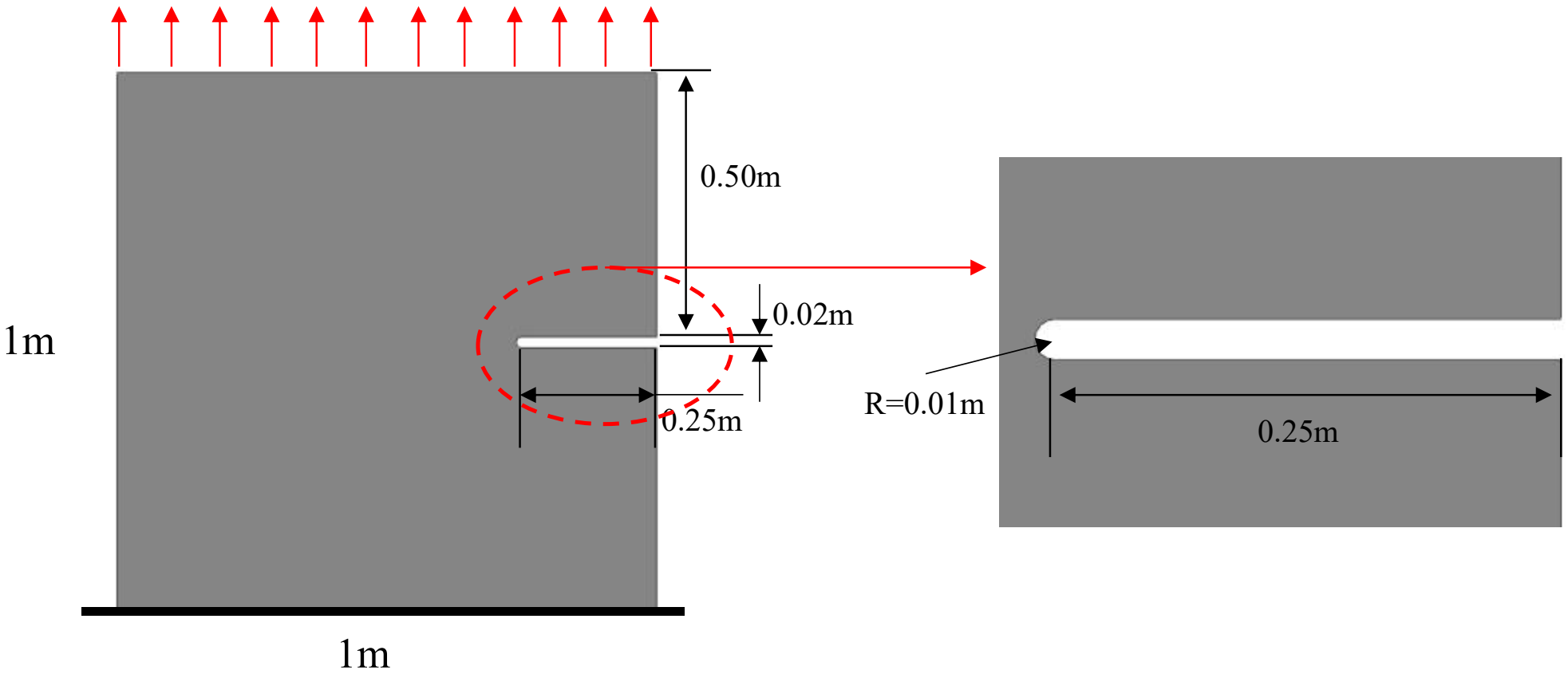

**Fig. 19. A plate with a long notch under tensile load: dimensions and boundary conditions.**

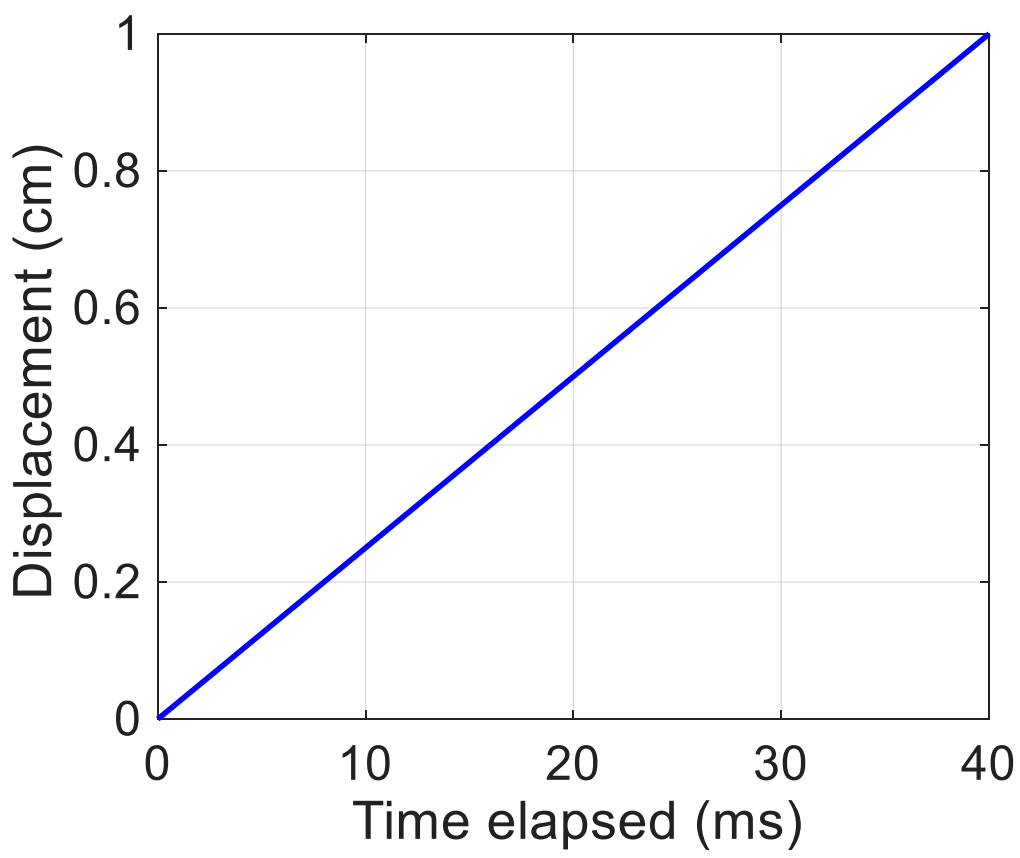

**Fig. 20. Load history for the plate with a notch.**

A node located at the notch tip is selected to monitor displacement and energy density. The results are presented in **Fig. 22**. As shown in **Fig. 22** (a), both the FEM and C-HiDeNN methods capture the displacement history with reasonable accuracy. However, C-HiDeNN demonstrates significantly higher accuracy in predicting strain energy density and von Mises stress at the tip, as illustrated in **Fig. 22** (b) and (c). The FEM results show increasing errors in both strain energy density and von Mises stress throughout the loading period, reaching maximum errors of 78.04% and approximately 50%, respectively, by the end of the loading step. In contrast, the C-HiDeNN method maintains an error of less than 5% for both quantities, indicating its superior predictive capability in regions of high stress concentration.

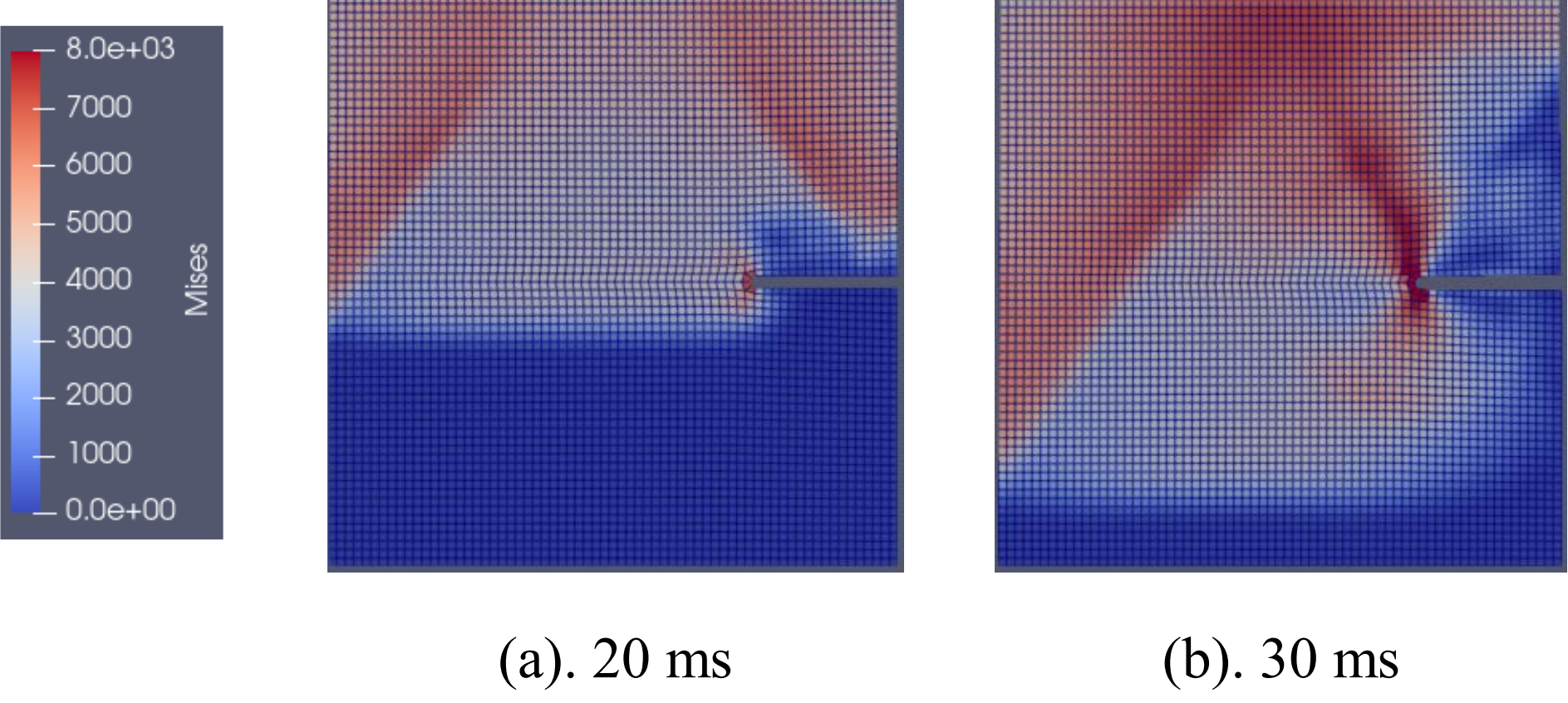

(a). 20 ms (b). 30 ms

**Fig. 21. Predicted von Mises Stress (Pa) distribution by C-HiDeNN at different time: (a) 20ms, (b) 30ms.**

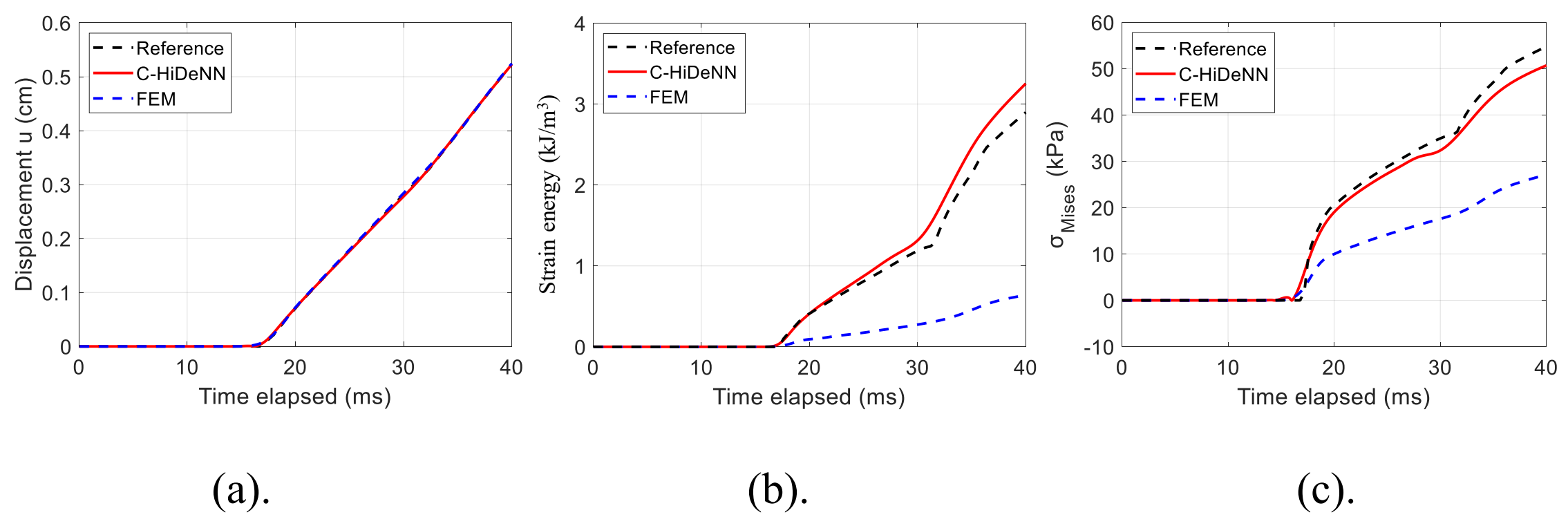

(a). (b). (c).

**Fig. 22. Comparison of results at the monitored point for the plate with an initial long notch problem among the reference solution, C-HiDeNN, and standard FEM: (a) Displacement *u* (cm), (b) Stran energy density (kJ/m³), (c) von Mises stress (kPa).**

### 4.5 3D hyperelastic Cook's beam with hole

In this example, we consider a 3D Cook's beam with a prescribed hole subjected to the compression load as shown in **Fig. 23** (a). The model dimensions are 0.5m by 0.7m by 0.1m as **Fig. 23** (a) shows, with a hole with radius R=0.1m located in the middle section. The material model is the same Neo-Hookean material model as in case 1 with parameters $C_{10}$ = 961.538 kPa, $D_1 = 4.8\times10^{-7}$ $\text{Pa}^{-1}$. The left surface is fixed, and a displacement load is applied to the right surface. The loading history is shown in **Fig. 23** (d), with a 1 mm displacement load applied over 40 ms. **Fig. 23** (b) shows the finite element mesh, consisting of 18,020 four-node tetrahedron elements (11,739 degrees of freedom), which is also used for generating the C-HiDeNN shape functions and performing the analysis. For comparison, **Fig. 23** (c) presents the reference solution mesh using the commercial code ABAQUS, which includes 19,520,836 C3D4 elements (four-node tetrahedron element), resulting in 10,391,217 degrees of freedom.

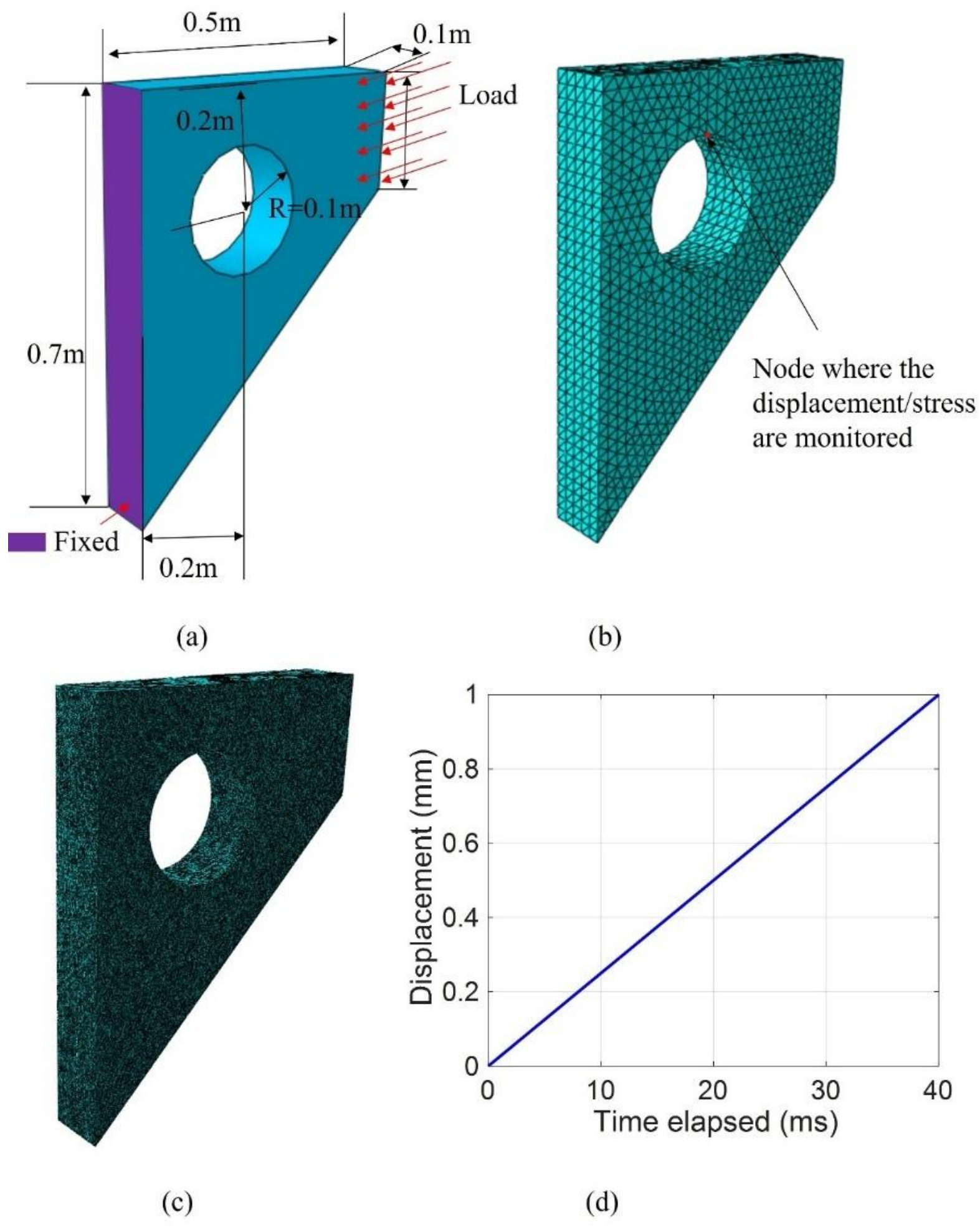

**Fig. 23. 3D hyperelastic beam with a hole: (a) Dimensions and boundary conditions (b) FE discretization (also for C-HiDeNN). (c) Reference solution mesh. (d) Load history.**

The results at the monitored point are shown in **Fig. 24**. Results show that C-HiDeNN method can capture displacement history (**Fig. 24** (a)) and the stress concentration (**Fig. 24** (b)) much better than finite element method with the same mesh. The C-HiDeNN approach maintains high accuracy while significantly reducing computational time. Reference solution needs around 61 hours, but C-HiDeNN method only takes 372 seconds, which leads to a speed-up of ~600.

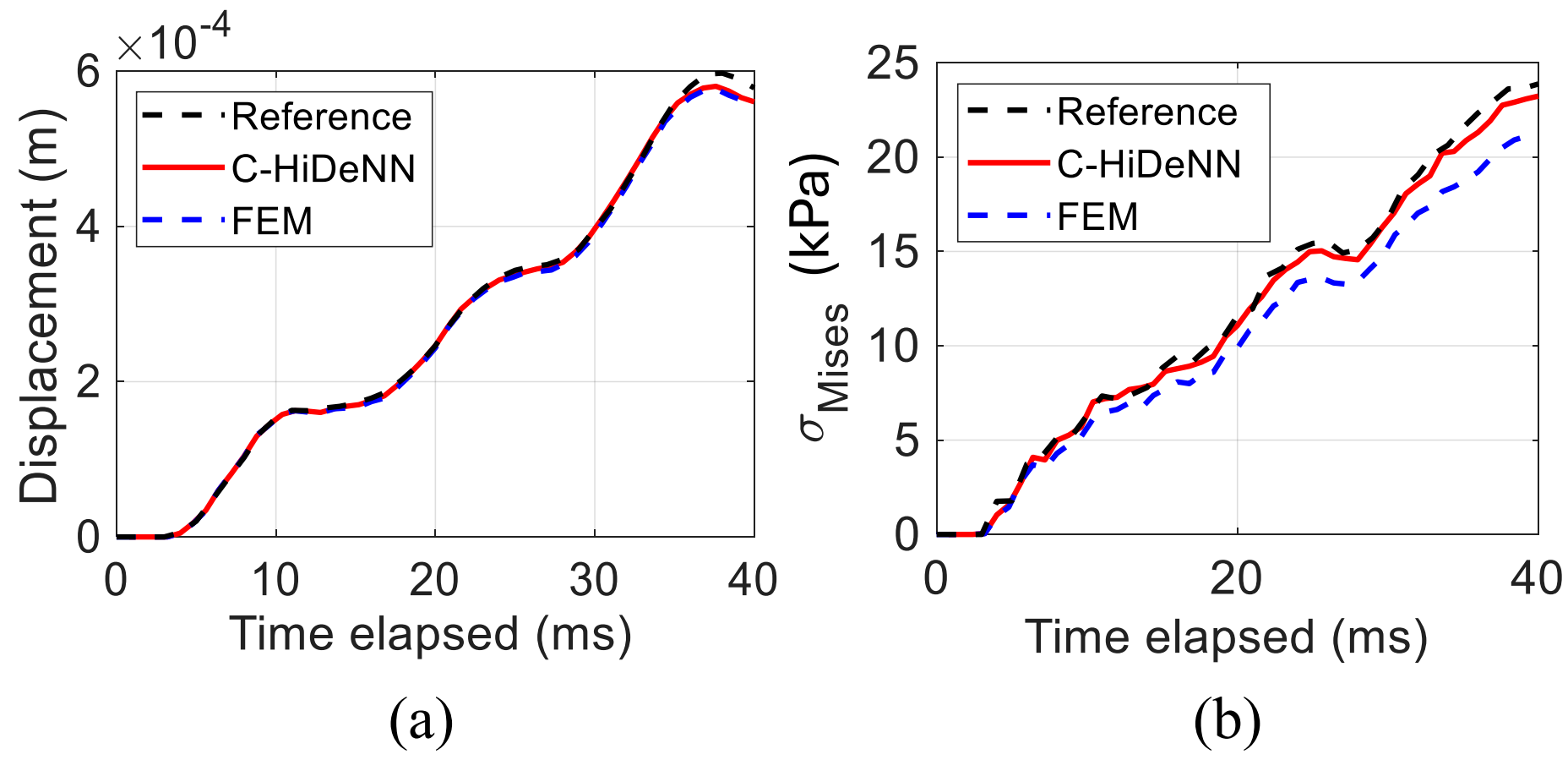

**Fig. 24. 3D rubber beam with a hole results comparison among reference solution, C-HiDeNN, and standard FEM: (a) displacement $u$ (m), (b) von Mises stress (kPa).**

### 4.6 3D mechanical bracket

In this example, we consider a 3D bracket subjected to the dynamic load as shown in **Fig. 25** (a). The model dimensions are shown as **Fig. 25** (a) with the fillet radius r=2.4mm. The material model is the same Neo-Hookean material model as case 1, with parameters $C_{10}$ = 961.538 kPa, $D_1 = 4.8\times10^{-7}$ Pa$^{-1}$. The bottom surface and four bolt holes are fixed, and displacement is applied to the bearing surface. The loading history is shown in **Fig. 25** (d), with a 0.02 mm displacement applied over 40 ms at a constant rate. The finite element mesh is shown in **Fig. 25** (b) and consists of 38,587 four-node tetrahedron elements (25,782 degrees of freedom). This mesh is also used for generating the C-HiDeNN shape functions and performing the analysis. For validation purposes, the reference solutions are obtained using a FE model generated in the commercial code ABAQUS with extremely refined mesh as shown in **Fig. 25** (c). This model has 20,895,415 C3D4 elements (four-node tetrahedron element), resulting in 11,026,611 degrees of freedom.

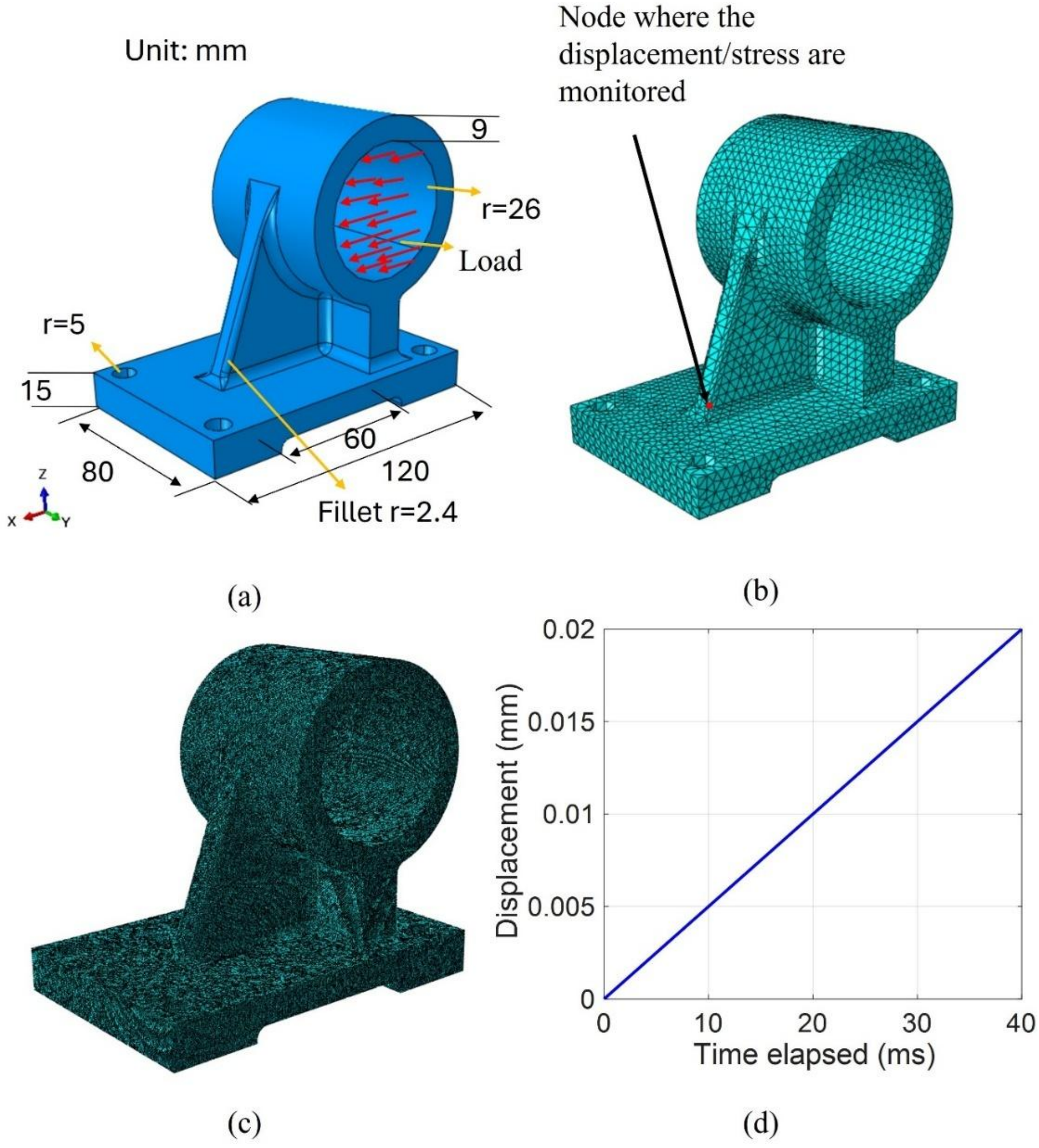

**Fig. 25. 3D simulation of a mechanical bracket under dynamic load: (a) Dimensions of the model, numbers shown are in mm unit. (b) FE discretization (also used for C-HiDeNN). (c) FE mesh for obtaining the reference solution. (d) load history.**

The results at the monitored point are shown in **Fig. 26**, which compares the predictions at a location exhibiting high stress concentration (as indicated in **Fig. 25** (b)) that demonstrates high stress concentration. The results show that C-HiDeNN accurately captures both displacement history (**Fig. 26** (a)) and stress concentration (**Fig. 26** (b)). Both predictions show significantly higher accuracy than those obtained using the finite element method (FEM) with the same mesh resolution. In addition to improved accuracy, C-HiDeNN offers substantial computational savings. While the reference solution requires approximately 43 hours to compute, C-HiDeNN completes the same task in just 14 minutes, resulting in an approximate 200× speed-up. A detailed comparison of the results at 30 ms with the reference solution is provided in **Table 2**. The data shows that C-HiDeNN achieves displacement and stress predictions with errors less than 5%, compared to errors exceeding 15% for FEM.

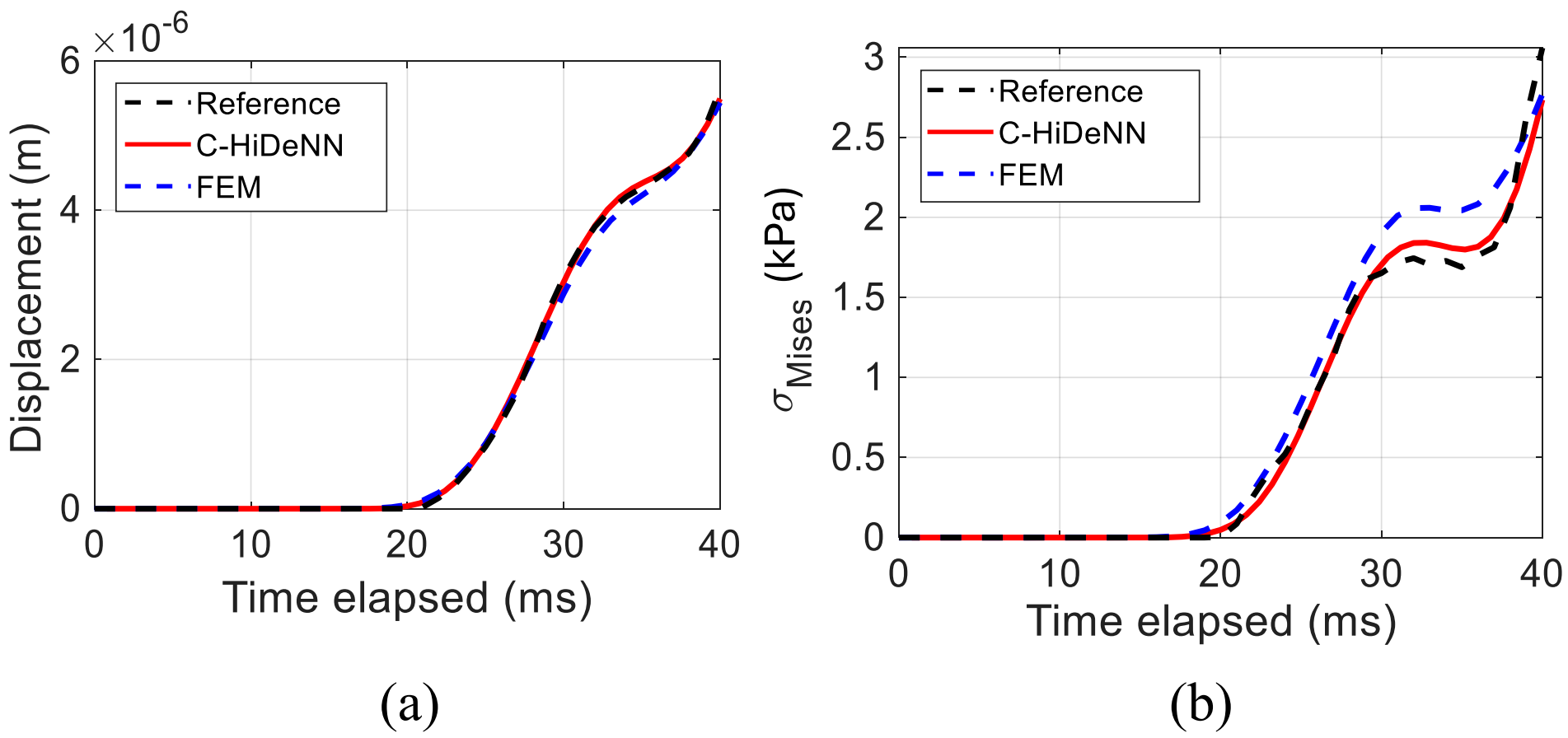

**Fig. 26. 3D mechanical bracket results comparison among reference solution, C-HiDeNN and standard FEM: (a) displacement $u$ (m), (b) von Mises stress (kPa).**

| Analysis/field result | FEM (C3D4) | C-HiDeNN |
|---|---|---|
| $u$ | 7.09% | 2.12% |
| $\boldsymbol{\sigma}_{mises}$ | 15.38% | 3.27% |

**Table 2: 3D mechanical bracket error analysis at 30ms (Reference solutions are obtained from fine mesh of ~20M C3D4 elements).**

## 5 Conclusion

In this work, we developed and demonstrated the application of Convolutional Hierarchical Deep-learning Neural Network Finite Element Method (C-HiDeNN) for nonlinear analysis, as a powerful extension of traditional finite element methods and the earlier HiDeNN framework. By incorporating convolution operators into shape function approximation, C-HiDeNN enables a hybrid FE approach that achieves higher-order accuracy, enhanced adaptivity, and significantly reduced computational costs.

Comprehensive numerical examples, including structural and biomedical applications under dynamic loading conditions, validated the accuracy and efficiency of C-HiDeNN. Across all cases, C-HiDeNN demonstrated superior accuracy in capturing displacement, stress concentration, and energy density compared to standard FEM, even when using significantly coarser meshes. In particular, the *s*-adaptivity feature enabled localized enrichment without modifying the global mesh, resulting in substantial reductions in computational time while maintaining high fidelity near critical regions such as notches or crack tips.

In conclusion, the C-HiDeNN framework offers a scalable, adaptive, and physically grounded approach to nonlinear dynamic analysis. Future work will focus on extending this method to integration with CAD representations, such as those based on tensor decomposition and NURBS, multi-physics coupling, and real-time digital twin applications.

## 6 Credit authorship contribution statement

Yingjian Liu: Writing-original draft, Software, Methodology, Formal analysis, Conceptualization, Writing – review & editing

Monish Yadav Pabbala: Software, Methodology, Writing – review & editing

Jiachen Guo: Methodology, Formal analysis, Conceptualization, Writing – review & editing

Chanwook Park: Methodology, Formal analysis, Conceptualization, Writing – review & editing

Gino Domel: Writing – review & editing, Conceptualization

Wing Kam Liu: Writing – review & editing, Project administration, Conceptualization

Dong Qian: Writing-original draft, Software, Methodology, Formal analysis, Conceptualization, Writing – review & editing, Project administration, Investigation.

## 7 Declaration of Competing Interest

The authors declare that they have no known competing financial interests or personal relationships that could have appeared to influence the work reported in this paper.